# MUTUAL INTERPRETABILITY OF WEAK ESSENTIALLY UNDECIDABLE THEORIES

Zlatan Damnjanovic

[Abstract: Kristiansen and Murwanashyaka recently proved that Robinson arithmetic, Q, is interpretable in an elementary theory of full binary trees,T. We prove that, conversely, T is interpretable in Q by producing a formal interpretation of T in an elementary concatenation theory $QT^+$, thereby also establishing mutual interpretability of T with several well-known weak essentially undecidable theories of numbers, strings and sets. We also introduce a "hybrid" elementary theory of strings and trees, WQT*, and establish its mutual interpretability with Robinson's weak arithmetic R, the weak theory of trees WT of Kristiansen and Murwanashyaka and the weak concatenation theory $WTC^{-\varepsilon}$ of Higuchi and Horihata.]

Key words: interpretability, full binary trees, Robinson arithmetic, concatenation theory, strings, essential undecidability

The classic monograph work of Tarski Mostowski and Robinson [8] isolated two weak formal theories of arithmetic, R and Q, as minimal "basis theories" for metamathematical arguments of foundational significance involving formalizing computation, incompleteness, undecidability, etc. The two theories were singled out as essentially undecidable, in that neither can consistently be extended to a decidable theory. The work introduced a powerful method for establishing incompleteness and undecidability of a wide range of mathematical theories built around the notion of relative interpretability of one theory in another. Roughly, a formula with a single free variable is chosen in the language of the second theory – the interpret<u>ing</u> theory -- to define the “universe of the interpretation”, and suitable definitions for the non-logical vocabulary of the first theory – the interpret<u>ed</u> theory -- are given in the language of the interpret<u>ing</u> theory. Formulae of the interpret<u>ed</u> theory are then translated into formulae of the interpret<u>ing</u> theory based on those definitions, in such a way that the logical operations are preserved under the translation and, crucially, all occurrences of quantifiers become relativized to the universe of the interpretation. Consequently, deductive relations between formulae are preserved: in particular, <u>theorems</u> of the interpret<u>ed</u> theory are translated into <u>theorems</u> of the interpret<u>ing</u> theory. In this specific sense reasoning in one theory is formally simulated in another theory, establishing relative consistency of the former in the latter. Once it is shown that R or Q is interpretable in some given theory, it follows from Tarski's methods that the latter is also essentially undecidable.

It was only within the last two decades that some light has been shed on what makes R and Q special, a result of work of many researchers, including (earlier work by) Collins and Halpern, Wilkie, Grzegorczyk, Zdanowski, Švejdar, Ganea, and, especially, Visser. One approach was to characterize them as mutually interpretable with concatenation theories (theories of strings) or weak subsystems of set theory, each naturally motivated and of independent interest in their own right (see [1] for further references). Another is to produce a "coordinate-free" characterization independent of a particular axiomatic presentation in some formal language, as, e.g., in the remarkable theorem of Visser [10]: a recursively axiomatizable theory is interpretable in R if and only if it is locally finitely satisfiable, that is, each finite subset of its non-logical axioms has a finite model.

An important new angle on these issues was recently introduced in the work of Kristiansen and Murwanashyaka [6]. They consider two elementary axiomatizations, WT and T, whose intended models are simple inductively generated structures like trees or terms, and rigorously develop a direct and novel approach to formalization of computation by ultra-elementary means. T is formulated in the language $\mathcal{L}_T = \{0, (\ ), \sqsubseteq\}$ with a single individual constant 0, a binary operation symbol ( , ) and a 2-place relational symbol $\sqsubseteq$ with the following axioms:

(T1) $\forall x,y\ \neg(x,y)=0,$

(T2) $\forall x,y,z,w\ [(x,y)=(z,w) \rightarrow x=z\ \&\ y=w]$

(T3) $\forall x\ [x\sqsubseteq 0 \leftrightarrow x=0]$

(T4) $\forall x,y,z\ [x\sqsubseteq(y,z) \leftrightarrow x=(y,z) \vee x\sqsubseteq y \vee x\sqsubseteq z]$

On the other hand, the theory WT is formulated in the same vocabulary, but has infinitely many axioms given by the two schemas

(WT1) $\neg(s=t)$ for any distinct variable-free terms s, t of $\mathcal{L}_T$,

(WT2) $\forall x\,(x \sqsubseteq t \leftrightarrow \bigvee_{s \in \mathcal{S}(t)} x=s)$ for each variable-free term t of $\mathcal{L}_T$, where $\mathcal{S}(t)$ is the set of all subterms of t.

The theory WT, which turns out to be contained in T, is proved to be mutually interpretable with R. The stronger theory T, which can be thought of as the basic theory of full binary trees, even though lacking induction is shown to be sufficiently strong to allow for a formal interpretation of basic arithmetical operations validating the axioms of Q. Kristiansen and Murwanashyaka further conjectured that, conversely, T is also formally interpretable in Q.

In this paper we prove that T is indeed interpretable in Q, by formally interpreting T in a theory of concatenation, $QT^+$, previously investigated in [1] and established to be mutually interpretable with Q along with a host of other theories whose intended interpretations are natural numbers, strings or sets. Hence T and Q are mutually interpretable. Further we formulate a weak theory of concatenation, $WQT^*$, and a "pseudo-concatenation" theory WQT, and establish their mutual interpretability with Robinson's R. (While R is deductively contained, hence also interpretable, in Q, the latter, being finitely axiomatized but having no finite model, by Visser's Theorem is not interpretable in R.)

Several distinct formulations of concatenation theory which have been put forward as standard axiomatizations and as such extensively studied are not deductively co-extensive. Some, like Grzegorczyk's theory TC, are centered around what came to be known as Tarski's Law (or Editor Axiom), and some of the variants include the empty string as a unit element. Others, such as the theory $QT^+$ used in [1] and here, and a closely related theory F originally introduced by Tarski in [8], are on their face more explicitly theories of semi-groups with two generators. Nonetheless, all these theories turn out to be mutually interpretable on account of their mutual interpretability with Q. Our choice of $QT^+$ is motivated by the "ground-up" approach exemplified in the formula-selection method expounded below in §3.

In §§1-2 we give a preview of our interpretation of T in concatenation theory. In §3 we introduce the concatenation theory $QT^+$, explain the main methodological tool used throughout the paper, the formula selection method applied to tractable strings and string forms, and develop elements of formal concatenation theory $QT^+$ related to tallies, adding of tallies and parts of strings. §4 we describe the essentials of the coding methods subsequently used in formalization of definitions by string recursion in §5. The resulting formal schema of definition is applied to obtain definitions of counting functions α and β which we rely on to construct the formal interpretation introduced in §§1-2. In §6 the interpretation is formally defined, and translations of the axioms of T formally verified. There we state the main result of the paper, the First Mutual Interpretability Theorem of Weak Essentially Undecidable Theories, relating T and $QT^+$ to a number of well-known theories of numbers, strings and sets. Finally, in §7 we introduce concatenation variants WQT and WQT* of Robinson's theory R and establish the corresponding Second Mutual Interpretability Theorem with the weak theory of trees WT.

Many of our arguments involve construction of specific formulas and tedious verifications of their specific properties. Most of these details can be found in the Appendix. The entire formal construction ultimately rests on coding of sets of strings by strings within $QT^+$, which is given in complete detail in [2]. We provide specific references as needed.

# 1. Trees as Strings

The intended domain of interpretation of the theory T is the set of variable-free $\mathcal{L}_T$-terms

(*) 0, (00), (0(00)), ((00)0), ((00)(00)),…

Alternatively, we may think of the domain as consisting of finite full binary trees – also called 2-trees -- trees in which every node other than the endnodes has two immediate descendants. In order to interpret T in concatenation theory, we need some way of representing these objects – terms or trees – by binary strings. We would like to do this directly, without having to rely on a coding of sets or sequences.

For this purpose we will use a variant of Polish notation to read binary strings as codes for inductively generated objects having the structure characteristic of terms or trees. Thus, e.g., the terms in (*) will be coded, respectively, by

(**) a, baa, babaa, bbaaa, bbaabaa, …

To obtain the string code from a given variable-free $\mathcal{L}_T$-term we proceed from left to right by replacing the left parentheses by b's and 0's by a's, ignoring the right parentheses.

Looking at the strings that are examples of term codes in (**), we note that they all share the following features:

(c1) the total number of a's in the string exceeds the total number of b's exactly by 1,

(c2) each proper initial segment of the string has at least as many b's as a's.

In other words, each of these strings is <u>its own smallest initial segment in which the number of *a*'s strictly exceeds the number of *b*'s.</u> We will take this to be the defining property of binary term/tree codes. We offer the following as informal justification. Each *b* indicates a branching vertex, incurring a "debt" of two "open places", which need to filled by completing the branchings. That can be done either immediately by simply writing *a*, an end node, or by opening another branching, temporarily increasing the "debt of open places". Each successive *a* reduces the "debt" of places to be filled by one, until all open branchings are completed and the last two remaining

"places" filled with *a*'s, resulting in a full binary tree. Ultimately, *b*'s in the binary code track the number of branchings, i.e. non-terminal nodes, and *a*'s the number of terminal nodes in the tree.

To define the domain of the formal interpretation of T in concatenation theory we will need to be able to single out by means of a formula of concatenation theory those among arbitrary strings that are term codes. Key role in this connection will be played by functions $\alpha$ and $\beta$ that count the number of occurrences of the letters a, b, resp., in a given binary string. They are defined as follows:

$\alpha(a) = 1$ $\quad\quad\quad\quad$ $\beta(a) = 0$

$\alpha(b) = 0$ $\quad\quad\quad\quad$ $\beta(b) = 1$

$\alpha(x*a) = \alpha(x)+1$ $\quad\quad\quad\quad$ $\beta(x*a) = \beta(x)$

$\alpha(x*b) = \alpha(x)$ $\quad\quad\quad\quad$ $\beta(x*b) = \beta(x)+1$

Call a string x is <u>almost even</u>, writing Æ(x), if (c1) $\alpha(x) = \beta(x)+1$, and (c2) for each proper initial segment u of x, $\alpha(u) \leq \beta(u)$.

Within concatenation theory the values of $\alpha$, $\beta$ will be expressed by b-tallies, i.e., strings of consecutive $b$'s. The functions $\alpha$ and $\beta$ are additive in that

$$\alpha(x*y) = \alpha(x)+\alpha(y) \quad \text{and} \quad \beta(x*y) = \beta(x)+\beta(y).$$

To express and verify these properties in concatenation theory we will need to introduce a suitable operation Addtally having the requisite properties of addition on non-negative integers. But the main problem to be solved is showing that $\alpha$ and $\beta$, which are defined by recursion on strings, can actually be defined in concatenation theory.

## 2. Outline of the Interpretation

The language $\mathcal{L}_C = \{$ a, b, *$\}$ of concatenation theory has two individual constants a, b, and a single binary operation symbol *. Its intended interpretation $\Sigma^*$ has as its domain the set of all non-empty finite strings of $a$'s and $b$'s, the constants 'a', 'b', resp., stand for the digits a, b (or 0, 1, resp.), and, for given strings x, y from the domain of $\Sigma^*$, we take x*y to be the string obtained by concatenation (i.e., juxtaposition) of the successive digits of y to the right of the end digit of x. Simply put, for variable-free terms s, t of $\mathcal{L}_C$, an atomic formula 's=t' is true in $\Sigma^*$ just in case s and t denote the very same binary string. For the purpose of informal exposition of the basic idea behind the interpretation we will avail ourselves, "as a first approximation", of formulations couched in the first-order theory $Th(\Sigma^*)$ consisting of all true sentences of $\mathcal{L}_C$ in $\Sigma^*$. Specifically, at this point we will simply <u>assume</u> that the graphs of the functions α, β, are expressible by some formulas $A^\#(x,y)$, $B^\#(x,y)$, resp., of $\mathcal{L}_C$ along with the graph of Addtally, and carry on reasoning informally within $Th(\Sigma^*)$. In subsequent sections we turn to the detailed technical work of actually proving these assumptions by formalizing string recursion in concatenation theory and verifying the corresponding translations into $\mathcal{L}_C$ of the axioms of T, all of which has to be formally carried out within an extremely weak subtheory $QT^+$ of $Th(\Sigma^*)$.

First, some abbreviations. Let $xBy \equiv \exists z\, x*z=y$ and $xEy \equiv \exists z\, z*x=y$.

Then let $x \subseteq_p y \equiv x=y \vee xBy \vee xEy \vee \exists y_1 \exists y_2\, y=y_1*(x*y_2)$.

(Often, we shall write xy for x*y.)

2.1(a) $\Sigma^* \vDash Æ(x) \rightarrow x=a \vee (bBx \,\&\, aaEx)$.

(b) $\Sigma^* \vDash Æ(x) \,\&\, x_2Ex \rightarrow \alpha(x_2) \geq \beta(x_2)+1$.

(c) $\Sigma^* \models Æ(x)\ \&\ Æ(u)\ \&\ xy=uv \rightarrow x=u\ \&\ u=v$.

Proof: (a) Clearly, $\Sigma^* \models Æ(a)$. Assume $\Sigma^* \models Æ(x)\ \&\ x \neq a$. Then $\Sigma^* \models \neg aBx$, by (c2). So $\Sigma^* \models bBx$. Note that $\Sigma^* \models \neg Æ(aa)\ \&\ \neg Æ(ab)\ \&\ \neg Æ(ba)\ \&\ \neg Æ(bb)$.

Hence any x such that $\Sigma^* \models Æ(x)$ must have a (proper) endsegment of length 2. Suppose $\Sigma^* \models x=x_1ab \vee x=x_1ba \vee x=x_1bb$, that is, $abEx \vee baEx \vee bbEx$.

By (c1) and (c2), $\Sigma^* \models \alpha(x) = \beta(x)+1$, and $\Sigma^* \models \alpha(x_1) \leq \beta(x_1)$. If $\Sigma^* \models abEx$ or $\Sigma^* \models baEx$, then $\Sigma^* \models \alpha(x) = \alpha(x_1)+1$ and $\Sigma^* \models \beta(x) = \beta(x_1)+1$. But then

$$\Sigma^* \models \alpha(x) = \beta(x)+1 = (\beta(x_1)+1)+1 = \beta(x_1)+2 \geq \alpha(x_1)+2 > \alpha(x_1)+1 = \alpha(x),$$

a contradiction. On the other hand, if $\Sigma^* \models bbEx$, then $\Sigma^* \models \alpha(x) = \alpha(x_1)$ and $\Sigma^* \models \beta(x) = \beta(x_1)+2$. But then

$$\Sigma^* \models \alpha(x) = \beta(x)+1 = (\beta(x_1)+2)+1 = \beta(x_1)+3 \geq \alpha(x_1)+3 = \alpha(x)+3 > \alpha(x),$$

a contradiction again. Hence $\Sigma^* \models \neg abEx\ \&\ \neg baEx\ \&\ \neg bbEx$. But then we must have $\Sigma^* \models aaEx$.

(b) Assume $\Sigma^* \models Æ(x)\ \&\ x_2Ex$. Then $\Sigma^* \models \exists x_1\ x=x_1x_2$, hence

$\Sigma^* \models \alpha(x_1) \leq \beta(x_1)$. But $\Sigma^* \models \alpha(x) = \beta(x)+1$ and

$$\Sigma^* \models \alpha(x) = \alpha(x_1x_2) = \alpha(x_1)+\alpha(x_2),$$

whereas $\Sigma^* \vDash \beta(x) = \beta(x_1x_2) = \beta(x_1)+\beta(x_2)$. Then

$$\Sigma^* \vDash \alpha(x_1)+\alpha(x_2) = \beta(x_1)+\beta(x_2)+1,$$

whence from $\Sigma^* \vDash \alpha(x_1) \leq \beta(x_1)$ we have $\alpha(x_2) \geq \beta(x_2)+1$, as claimed.

(c) Assume $\Sigma^* \vDash Æ(x) \ \&\ Æ(u) \ \&\ xy=uv$. We have that

$$\Sigma^* \vDash (x=u \ \&\ y=v) \vee xBu \vee uBx.$$

Suppose $\Sigma^* \vDash xBu$. From $\Sigma^* \vDash Æ(u)$, $\Sigma^* \vDash \alpha(x) \leq \beta(x)$, and from $\Sigma^* \vDash Æ(x)$, $\Sigma^* \vDash \alpha(x) = \beta(x)+1$. But then $\Sigma^* \vDash \beta(x)+1 \leq \beta(x)$, a contradiction.

Likewise if $\Sigma^* \vDash uBx$. Hence $\Sigma^* \vDash x=u \ \&\ y=v$. ■

2.2 $\Sigma^* \vDash Æ(x) \leftrightarrow x=a \vee \exists! y,z\ (x=b(yz) \ \&\ Æ(y) \ \&\ Æ(z))$.

<u>Proof</u>: (⇐) Assume $\Sigma^* \vDash Æ(y) \ \&\ Æ(z) \ \&\ x=byz$. Then

$$\Sigma^* \vDash \alpha(y) = \beta(y)+1 \ \&\ \alpha(z) = \beta(z)+1.$$

Now, $\Sigma^* \vDash \alpha(x) = \alpha(byz) = \alpha(yz) = \alpha(y)+\alpha(z)$

and $\Sigma^* \vDash \beta(x) = \beta(byz) = \beta(b)+\beta(yz) = \beta(y)+\beta(z)+1$. Then

$$\Sigma^* \vDash \alpha(x) = \alpha(y)+\alpha(z) = (\beta(y)+1)+(\beta(z)+1) = (\beta(y)+\beta(z)+1)+1 = \beta(x)+1$$

which verifies (c1). For (c2), assume $\Sigma^* \models uBx$, i.e., $\Sigma^* \models uBbyz$.

Then $\Sigma^* \models u{=}b \ \vee \ uBby \ \vee u{=}by \ \vee \exists z_1(z_1Bz \ \& \ u{=}byz_1)$.

To illustrate the proof, we consider the case $\Sigma^* \models \exists z_1(z_1By \ \& \ u{=}byz_1)$.

Then from $\Sigma^* \models Æ(z)$, $\Sigma^* \models \alpha(z_1) \leq \beta(z_1)$, and from $\Sigma^* \models Æ(y)$,

$\Sigma^* \models \alpha(y) = \beta(y)+1$. Then $\Sigma^* \models \alpha(u) = \alpha(byz_1) = \alpha(yz_1) = \alpha(y)+\alpha(z_1)$ and

$\Sigma^* \models \ \beta(u) = \beta(byz_1) = \beta(b)+\beta(yz_1) = \beta(y)+\beta(z_1)+1$. Hence

$\Sigma^* \models \alpha(u) = \alpha(y)+\alpha(z_1) = (\beta(y)+1)+\alpha(z_1) \leq$

$$\leq (\beta(y)+1)+\beta(z_1) = \beta(y)+\beta(z_1)+1 = \beta(u).$$

Thus $\Sigma^* \models \ \alpha(u) \leq \beta(u)$. This completes the proof of (c2). So $\Sigma^* \models Æ(x)$.

($\Rightarrow$) Assume $\Sigma^* \models Æ(x) \ \& \ x{\neq}a$. Then, by 2.1(a), $\Sigma^* \models bBx \ \& \ aaEx$, that is,

$$\Sigma^* \models \exists x_1 \, x{=}bx_1 \ \& \ \exists x_2 \, x{=}x_2 \, aa.$$

So $\Sigma^* \models bx_1{=}x_2aa$. We may assume that $\Sigma^* \models bBx_2$, for if $\Sigma^* \models x_2{=}b$, then $\Sigma^* \models \ x{=}b(aa)$ and we may take $y{=}a$ and $z{=}a$. So $\Sigma^* \models \exists x_3 \, x_2{=}bx_3$, and $\Sigma^* \models x{=}bx_1{=}x_2(aa){=}bx_3(aa)$, whence $\Sigma^* \models x_1{=}x_3(aa)$. Let $y_j$ be a proper initial segment of $x_1$, and $z_j$ the corresponding endsegment of $x_1$ such that $\Sigma^* \models \ y_jz_j{=}x_1$. At least one $y_j$ has the property

(*) $\qquad \Sigma^* \models \alpha(y_j) = \beta(y_j)+1.$

Consider, e.g., $x_3a$. From hypothesis $\Sigma^* \models Æ(x)$ we have $\Sigma^* \models \alpha(x) = \beta(x)+1$.

But $\Sigma^* \vDash \alpha(x) = \alpha(b((x_3a)a)) = \alpha(b)+\alpha(x_3a)+\alpha(a) = \alpha(x_3a)+1$ and

$\Sigma^* \vDash \beta(x) = \beta(b((x_3a)a)) = \beta(b)+\beta(x_3a)+\beta(a) = 1+\beta(x_3a)$.

Then $\Sigma^* \vDash \alpha(x_3a) = \alpha(x)-1 = \beta(x) = \beta(x_3a)+1$.

Let $y_i$ be the <u>shortest</u> initial segment of $x_1$ with the property (*). Then

$$\Sigma^* \vDash x_1 = y_iz_i \ \& \ \alpha(y_i) = \beta(y_i)+1.$$

We claim that (i) $\Sigma^* \vDash \alpha(z_i) = \beta(z_i)+1$, (ii) $\Sigma^* \vDash \forall u\,(uBy_i \rightarrow \alpha(u)\leq\beta(u))$,

and (iii) $\Sigma^* \vDash \forall v\,(vBz_i \rightarrow \alpha(v)\leq\beta(v))$.

For (i) we have $\Sigma^* \vDash \alpha(x) = \alpha(bx_1) = \alpha(x_1) = \alpha(y_iz_i) = \alpha(y_i)+\alpha(z_i)$

and $\Sigma^* \vDash \beta(x) = \beta(bx_1) = 1+\beta(x_1) = 1+\beta(y_iz_i) = 1+ \beta(y_i)+\beta(z_i)$.

Then $\Sigma^* \vDash \alpha(y_i)+\alpha(z_i) = (1+ \beta(y_i)+\beta(z_i))+1$, and from $\Sigma^* \vDash \alpha(y_i) = \beta(y_i)+1$

we obtain $\Sigma^* \vDash \alpha(z_i) = \beta(z_i)+1$.

For (ii), suppose $\Sigma^* \vDash uBy_i$. Since $\Sigma^* \vDash x_1 = y_iz_i$, we then have $\Sigma^* \vDash uBx_1$. But then, by the choice of $y_i$, $\Sigma^* \vDash \alpha(u)\leq\beta(u)$. For (iii), suppose $\Sigma^* \vDash vBz_i$. Then $\Sigma^* \vDash \exists w\ z_i = vw$, whence $\Sigma^* \vDash wEx$. From $\Sigma^* \vDash Æ(x)$, by 2.1(b),

$\Sigma^* \vDash \alpha(w) \geq \beta(w)+1$. But $\Sigma^* \vDash \alpha(z_i) = \alpha(v)+\alpha(w)$ and

$\Sigma^* \vDash \beta(z_i) = \beta(v)+\beta(w)$. By (i), $\Sigma^* \vDash \alpha(v)+\alpha(w) = \beta(v)+\beta(w)+1$.

Then from $\Sigma^* \vDash \alpha(w) \geq \beta(w)+1$, we have $\Sigma^* \vDash \alpha(v)\leq\beta(v)$.

From (i)-(iii) we have that $\Sigma^* \vDash Æ(y_i)$ & $Æ(z_i)$. The uniqueness of y, z follows from 2.1(c).■

The proof of 2.2 yields an algorithm for extracting the description of a tree from a given Æ string *x*: (i) Drop the initial *b*. (ii) If the next digit is *a*, that is the left node Æ string; the rest of the string is the right node Æ string. (iii) If the next digit is *b*, take the shortest initial segment *y* of the remainder of the original string such that α(*y*)=β(*y*)+2; then the string *by* is the left node Æ string, and the endsegment of the remainder corresponding to *by* is the right node Æ string. Repeat steps (i)-(iii) until no *b*'s are left.

2.3 $\Sigma^* \vDash Æ(x)$ & $Æ(y)$ & $Æ(z) \rightarrow (x \subseteq_p byz \rightarrow x=byz \text{ v } x \subseteq_p y \text{ v } x \subseteq_p z)$.

<u>Proof</u>: Assume $\Sigma^* \vDash x \subseteq_p byz$ where $\Sigma^* \vDash Æ(x)$ & $Æ(y)$ & $Æ(z)$. Now, we have that $\Sigma^* \vDash x=byz \text{ v } x=b \text{ v } x \subseteq_p yz \text{ v } \exists u(uByz \;\&\; x=bu)$.

Suppose that $\Sigma^* \vDash \exists u(uByz \;\&\; x=bu)$. From $\Sigma^* \vDash Æ(y)$ & $Æ(z)$, by 2.2, $\Sigma^* \vDash Æ(byz)$. From $\Sigma^* \vDash uByz$, we have $\Sigma^* \vDash \exists v\, uv=yz$, whence $\Sigma^* \vDash buBb(yz)$. Thus $\Sigma^* \vDash xBb(yz)$. But from $\Sigma^* \vDash Æ(byz)$, $\Sigma^* \vDash \alpha(x) \leq \beta(x)$, which contradicts $\Sigma^* \vDash Æ(x)$. So $\Sigma^* \vDash \exists u(uByz \;\&\; x=bu)$ is ruled out.
By 2.1(a), so is $\Sigma^* \vDash x=b$. So we are left with $\Sigma^* \vDash x \subseteq_p byz \rightarrow x=byz \text{ v } x \subseteq_p yz$.
Supposing $\Sigma^* \vDash x \subseteq_p yz$, we have that $\Sigma^* \vDash x=yz \ \text{v} \ x\subseteq_p y \ \text{v} \ x\subseteq_p z \ \text{v} \ \exists y_1(y_1Ey \ \& \ x=y_1z) \ \text{v}$

$$\text{v} \ \exists z_1 (z_1Bz \ \& \ x=yz_1) \ \text{v} \ \exists y_1,z_1 (y_1Ey \ \& \ z_1Bz \ \& \ x=y_1z_1).$$

Assume $\Sigma^* \vDash x=yz$. Then from $\Sigma^* \vDash Æ(y) \ \& \ Æ(z)$, we have $\Sigma^* \vDash \alpha(y) = \beta(y)+1$ and $\alpha(z) = \beta(z)+1$. But $\Sigma^* \vDash \alpha(yz) = \alpha(y)+\alpha(z)$, so

$$\Sigma^* \vDash \alpha(yz) = (\beta(y)+1)+(\beta(z)+1) = \beta(y)+\beta(z)+2.$$

On the other hand, $\Sigma^* \vDash \beta(yz) = \beta(y)+\beta(z)$. Thus $\Sigma^* \vDash \alpha(yz) = \beta(yz)+2$, whence from $\Sigma^* \vDash x=yz$, we derive $\Sigma^* \vDash \alpha(x) = \beta(x)+2$, contradicting $\Sigma^* \vDash Æ(x)$. So $\Sigma^* \vDash x=yz$ is ruled out.

Suppose now that $\Sigma^* \vDash \exists y_1(y_1Ey \ \& \ x=y_1z)$, so $\Sigma^* \vDash y_1Bx$. From $\Sigma^* \vDash Æ(x)$, $\Sigma^* \vDash \alpha(y_1) \leq \beta(y_1)$. But from $\Sigma^* \vDash Æ(y) \ \& \ y_1Ey$, we obtain, by 2.1(b), $\Sigma^* \vDash \alpha(y_1) \geq \beta(y_1)+1$, a contradiction.

Suppose that $\Sigma^* \vDash \exists z_1 (z_1Bz \ \& \ x=yz_1)$, so $\Sigma^* \vDash yBx$. But then from $\Sigma^* \vDash Æ(x)$, we have $\Sigma^* \vDash \alpha(y) \leq \beta(y)$, and from $\Sigma^* \vDash Æ(y)$, $\Sigma^* \vDash \alpha(y) = \beta(y)+1$, again a contradiction. If $\Sigma^* \vDash \exists y_1,z_1 (y_1Ey \ \& \ z_1Bz \ \& \ x=y_1z_1)$, we derive a contradiction by reasoning as in either of the two preceding cases.

The other cases having been ruled out, we conclude under the principal hypothesis that $\Sigma^* \vDash x\subseteq_p yz \rightarrow x\subseteq_p y \ \text{v} \ x\subseteq_p z$, and further that

$$\Sigma^* \vDash x\subseteq_p byz \rightarrow x=byz \ \text{v} \ x\subseteq_p y \ \text{v} \ x\subseteq_p z,$$

as required.■

If we take the domain to consists of Æ strings, 2.1(c), 2.2 and 2.3 suffice to give the “first approximation” of our interpretation of T in concatenation theory: translations of (T1)-(T4) will be validated in $\Sigma^*$ if we model the term/tree-building operation x, y ↦ (xy) by *b*xy, the subterm/subtree relation ⊑ by the substring relation $\subseteq_p$ between Æ strings, and the digit *a* is taken to stand for the simple term 0. The entire project, however, hinges on definability of the counting functions α and β in concatenation theory. Showing that the latter contains resources needed to formally justify definitions by elementary recursion on strings requires, first, that we precisely formulate concatenation theory as a formal theory, and second, that we introduce codings for ordered pairs of strings, sequences of such, etc., and verify their properties relevant to the argument in that formal theory. We now turn to that task. In the process we shall make crucial use of the method of formula selection explained in [1].

## 3. Formal Concatenation Theory

We shall work within a first-order theory formulated in $\mathcal{L}_C = \{ a, b, * \}$, with the universal closures of the following conditions as axioms:

(QT1) x*(y*z)=(x*y)*z

(QT2) $\neg(x*y=a)$ & $\neg(x*y=b)$

(QT3) $(x*a=y*a \rightarrow x=y)$ & $(x*b=y*b \rightarrow x=y)$ &

& $(a*x=a*y \rightarrow x=y)$ & $(b*x=b*y \rightarrow x=y)$

(QT4) $\neg(a*x=b*y)$ & $\neg(x*a=y*b)$

(QT5) $x=a \vee x=b \vee (\exists y(a*y=x \vee b*y=x) \;\&\; \exists z(z*a=x \vee z*b=x))$

On account of (QT1), we sometimes omit parentheses and * when writing

(x*y).

It is convenient to have a function symbol for a successor operation on strings:

(QT6) $Sx=y \leftrightarrow ((x=a \;\&\; y=b) \vee (\neg x=a \;\&\; x*b=y))$.

Since (QT6) is basically a definition, adding it to the rest results in an inessential (i.e. conservative) extension. We call this theory $QT^+$.

Let $xRy \equiv (x=a \;\&\; \neg y=a) \vee xBy$.

Provably in $QT^+$, $xRy \vee x=y$ is a discrete preordering of strings (see [1]).

We shall call a formula I(x) in the language of $QT^+$ a <u>string form</u> if

$QT^+ \vdash I(a)$, $QT^+ \vdash I(b)$, $QT^+ \vdash I(x) \rightarrow I(x*a)$ and $QT^+ \vdash I(x) \rightarrow I(x*b)$.

(Note: in [1] and [2] such formulae were called string concepts.) String forms will allow us to restrict our attention, systematically step-by-step, to strings that satisfy conditions expressible by specifically selected formulas provided the latter can be proved in $QT^+$ to apply to "sufficiently many" strings. We say that a string form J is <u>stronger</u> <u>than</u> I if $QT^+ \vdash \forall x\,(J(x) \rightarrow I(x))$ and write $J \subseteq I$.

Let $I_0(x) \equiv \forall y\,(yRx \vee y=x \rightarrow \neg yRy)$. We call $I_0$ strings <u>tractable</u>.

3.1(a) $I_0(x)$ is a string form.

(b) For any string form $I \subseteq I_0$ there is a string form $J \subseteq I$ such that

$$QT^+ \vdash \forall x\,\forall y\,(J(x)\;\&\;J(y) \rightarrow J(x*y)).$$

(c) For any string form $I \subseteq I_0$ there is a string form $J_{\leq} \subseteq I$ such that

$$QT^+ \vdash \forall x\,(J_{\leq}(x)\;\&\;y \leq x \rightarrow J_{\leq}(y)).$$

(d) For any string form $I \subseteq I_0$ there is a string form $J \subseteq I$ such that

$$QT^+ \vdash \forall x \in J\,\forall y\,(y \subseteq_p x \rightarrow J(y)).$$

(e) For any string form $I \subseteq I_0$ there is a string form $J \equiv I_{LC} \subseteq I$ such that

$$QT^+ \vdash \forall z \in J\,\forall x,y\,(z*x = z*y \rightarrow x=y).$$

(f) For any string form $I \subseteq I_0$ there is a string form $J \subseteq I$ such that

$$QT^+ \vdash \forall z \in J\,\forall x,y\,(x*z = y*z \rightarrow x=y).$$

For proofs, see [1], and [2], (3.2), (3.3), (3.13), (3.7) and (3.6).■

Parts (b)-(c) tell us that when establishing that a given string form I may be strengthened to a string form J with another property, we can always strengthen the string form J to one that is also closed with respect to * or downward closed with respect to $\leq$ or $\subseteq_p$ .

We define $\quad \mathrm{Tally}_a(x) \equiv \forall y \subseteq_p x\ (\mathrm{Digit}(y) \rightarrow y=a)$

and $\quad \mathrm{Tally}_b(x) \equiv \forall y \subseteq_p x\ (\mathrm{Digit}(y) \rightarrow y=b)$ where $\quad \mathrm{Digit}(x) \equiv x=a \vee x=b$.

Write $x<y$ for $I_0(x)\ \&\ I_0(y)\ \&\ xRy$. As usual, $x\leq y$ stands for $x<y \vee x=y$.

The following properties of tallies are easily established:

3.2 (a) $QT^+ \vdash \mathrm{Tally}_b(y) \rightarrow \mathrm{Tally}_b(Sy)$.

(b) $QT^+ \vdash \mathrm{Tally}_b(y) \leftrightarrow y=b \vee \exists y_1\ (\mathrm{Tally}_b(y_1)\ \&\ y=Sy_1)$.

(c) $QT^+ \vdash \forall v,u\ (\mathrm{Tally}_b(v)\ \&\ u<v \rightarrow Su\leq v)$.

(d) $QT^+ \vdash \mathrm{Tally}_b(y) \rightarrow (x<y \leftrightarrow Sx<Sy)$.

For some further properties we have to resort to string forms:

3.3(a) For any string form $I\subseteq I_0$ there is a string form $J\equiv I_{CTC}\subseteq I$ such that

$QT^+ \vdash \forall z \in J\ \forall y\ (\mathrm{Tally}_b(y)\ \&\ \mathrm{Tally}_b(z) \rightarrow \mathrm{Tally}_b(y*z))$.

(b) For any string form $I\subseteq I_0$ there is a string form $J\subseteq I$ such that

$QT^+ \vdash \forall z \in J\ \forall x\ (Tally_b(x)\ \&\ Tally_b(z)\ \rightarrow\ x \leq z\ \ v\ \ z \leq x)$.

(c) For any string form $I \subseteq I_0$ there is a string form $J \equiv I_{3.3(c)} \subseteq I$ such that

$$QT^+ \vdash \forall u \in J\ (Tally_b(u) \rightarrow\ u{*}b = b{*}u).$$

(d) For any string form $I \subseteq I_0$ there is a string form $J \subseteq I$ such that

$$QT^+ \vdash \forall y \in J\ \forall x\ (Tally_b(x)\ \&\ Tally_b(y)\ \rightarrow\ Sx{*}y = x{*}Sy = S(x{*}y)).$$

(e) For any string form $I \subseteq I_{3.3(c)}$ there is a string form $J \equiv I_{COMM} \subseteq I$ such that

$$QT^+ \vdash \forall u,v \in J\ (Tally_b(u)\ \&\ Tally_b(v)\ \rightarrow\ u{*}v = v{*}u).$$

For proofs, see [2], (4.5), (4.6), (4.8) and (4.10).■

Let Addtally(x,y,z) abbreviate the formula

$$(Tally_b(x)\ \&\ Tally_b(y)\ \&\ ((x=b\ \&\ z=y)\ v\ (y=b\ \&\ z=x)\ v$$
$$v\ \exists x_1,y_1(Tally_b(x_1)\ \&\ x=Sx_1\ \&\ Tally_b(y_1)\&\ y=Sy_1\ \&\ z=x{*}y_1))\ v$$
$$v\ ((\neg Tally_b(x)\ v\ \neg Tally_b(y))\ \&\ z=b)$$

We want to show that, provably in $QT^+$, Addtally(x,y,z) behaves like the graph of addition function on natural numbers. The following are immediate consequences of definitions:

3.4(a) $QT^+ \vdash$ Addtally(x,y,v) & Addtally(x,y,w) $\rightarrow$ v=w.

(b) $QT^+ \vdash Tally_b(x) \rightarrow$ Addtally(x,b,x). ("x+0 = x")

(c) $QT^+ \vdash Tally_b(y) \rightarrow$ Addtally(b,y,y). ("0+y = y")

(d) $QT^+ \vdash Tally_b(x) \to Addtally(x,bb,Sx)$. ("x+1 = Sx")

(e) $QT^+ \vdash Tally_b(x) \,\&\, Tally_b(y) \to (Addtally(x,y,z) \to Addtally(x,yb,zb))$.

("x+Sy = S(x+y)")

We also have:

3.5(a) For any string form $I \subseteq I_0$ there is a string form $J \equiv I_{Add} \subseteq I$ such that

$$QT^+ \vdash \forall x,y \in J\ \exists! z \in J\ (Tally_b(z) \,\&\, Addtally(x,y,z)).$$

(b) $QT^+ \vdash \forall z \in I_0\ (Tally_b(u) \,\&\, Tally_b(v) \,\&$

$$\&\, Addtally(x,u,y) \,\&\, Addtally(x,v,z) \,\&\, u \leq v \to y \leq z).$$

("$u \leq v \to x+u \leq x+v$")

(c) For any string form $I \subseteq I_0$ there is a string form $J \subseteq I$ such that

$$QT^+ \vdash \forall y \in J\ (Tally_b(y) \to Addtally(bb,y,Sy).$$ ("1+y = Sy")

(d) For any string form $I \subseteq I_0$ there is a string form $J \subseteq I$ such that

$$QT^+ \vdash \forall y \in J\ \forall x,z\ (Tally_b(x) \,\&\, Tally_b(y) \,\&\, Addtally(x,y,z) \to Addtally(xb,y,zb))$$

("Sx+y = S(x+y")

(e) For any string form $I \subseteq I_0$ there is a string form $J \subseteq I$ such that

$$QT^+ \vdash \forall x \in J\ \forall y,z,v\ (Tally_b(x) \,\&\, Tally_b(y) \,\&\, Tally_b(z) \to$$

$$\to (Addtally(x,y,v) \,\&\, Addtally(x,z,v) \to y=z)).$$

("$x+y=x+z \to y=z$")

(f) For any string form $I \subseteq I_0$ there is a string form $J \subseteq I$ such that

$$QT^+ \vdash \forall y \in J\ \forall x\ (Tally_b(x) \,\&\, Tally_b(y) \to$$

$$\to (x \leq y \leftrightarrow \exists z (Tally_b(z) \,\&\, Addtally(z,x,y)))).$$

("$x \leq y \leftrightarrow \exists z\, z+x=y$")

(g) For any string form $I \subseteq I_0$ there is a string form $J \subseteq I$ such that

$$QT^+ \vdash \forall x,y \in J\ (\mathrm{Addtally}(x,y,z) \rightarrow \mathrm{Addtally}(y,x,z)). \qquad (\text{"}x+y=y+x\text{"})$$

(h) For any string form $I \subseteq I_0$ there is a string form $J \subseteq I$ such that

$$QT^+ \vdash \forall x,y,z \in J\ (\mathrm{Addtally}(x,y,u)\ \&\ \mathrm{Addtally}(u,z,v_1)\ \&\ \mathrm{Addtally}(y,z,w)\ \&$$
$$\&\ \mathrm{Addtally}(x,w,v_2) \rightarrow v_1=v_2)$$

("$(x+y)+z=x+(y+z)$")

(i) For any string form $I \subseteq I_0$ there is a string form $J \subseteq I$ such that

$$QT^+ \vdash \forall x_2,y_1,y_2 \in J\ \forall x_1,z_1,z_2\ (\mathrm{Tally}_b(x_2)\ \&\ \mathrm{Tally}_b(y_1)\ \&\ \mathrm{Tally}_b(y_2)\ \&$$
$$\&\ \mathrm{Addtally}(x_1,x_2,z_1)\ \&\ \mathrm{Addtally}(y_1,y_2,z_2)\ \&\ x_1 \leq y_1\ \&\ z_1=Sz_2 \rightarrow$$
$$\rightarrow Sy_2 \leq x_2).$$

("$x_1+x_2 = (y_1+y_2)+1\ \&\ x_1 \leq y_1 \rightarrow y_2+1 \leq x_2$")

Proof: For (a), let $J \equiv I_{CTC}$ from 3.3(a). For (c) and (d), let J be as in 3.3(c). For (e), let $J \equiv I_{LC}$ from 3.1(d). For (f) and (g), let $J \equiv I_{COMM}$ from 3.3(e). For (h), let $J \equiv J_1\ \&\ J_2$ where $J_1$ is $I_{CTC}$ and $J_2$ as in 3.3(c). Finally, for (i), let $J \equiv I_{LC}\ \&\ I_{CTC}\ \&\ I_{3.3(c)}\ \&\ I_{COMM}$ and see Appendix.■

We now turn to the part-of relation $\subseteq_p$ between strings. To prevent unpleasant surprises, we want to make sure that this relation has natural properties we would normally expect it to have.

3.6(a) $QT^+ \vdash x \subseteq_p y \ \& \ y \subseteq_p z \rightarrow x \subseteq_p z$.

(b) For any string form $I \subseteq I_0$ there is a string form $J \subseteq I$ such that

$$QT^+ \vdash \forall x \in J \ \neg xEx.$$

(c) For any string form $I \subseteq I_0$ there is a string form $J \subseteq I$ such that

$$QT^+ \vdash \forall x \in J \ \neg \exists x_1, x_2 \ (x_1 x x_2 = x).$$

(d) For any string form $I \subseteq I_0$ there is a string form $J \subseteq I$ such that

$$QT^+ \vdash \forall x \in J \ \forall y \ (x \subseteq_p y \ \& \ y \subseteq_p x \ \rightarrow \ x = y).$$

(e) For any string form $I \subseteq I_0$ there is a string form $J \subseteq I$ such that

$$QT^+ \vdash \forall x \in J \ \forall y \ (\neg xy \subseteq_p x \ \& \ \neg yx \subseteq_p x).$$

Proof: For (b) and (c), see [2], (3.4) and (3.5). For (d) and (e), see [2], (3.11) and (3.12).■

We now specifically consider proper initial segments and endsegments. The initial segments of arbitrary strings can be totally ordered by the initial-segment-of relation B, rendering the partial ordering < in which *a* is the least element tree-like:

3.7(a) For any string form $I \subseteq I_0$ there is a string form $J_{LOIS} \subseteq I$ such that

$$QT^+ \vdash \forall x \in J \ \forall u, v \ (uBx \ \& \ vBx \ \rightarrow \ u = v \ \vee \ uBv \ \vee \ vBu).$$

(b) For any string form $I \subseteq I_0$ there is a string form $J \subseteq I$ such that

$$QT^+ \vdash \forall y,z \in J\ \forall x\ (xByz \leftrightarrow xBy \lor x=y \lor \exists w(wBz \ \&\ yw=x)).$$

(c) For any string form $I \subseteq I_0$ there is a string form $J \subseteq I$ such that

$$QT^+ \vdash \forall x,y \in J\ \forall u\ (uBb(xy) \rightarrow u=b \lor uBbx \lor u=bx \lor \exists y_1(y_1By \ \&\ u=bxy_1)).$$

(d) For any string form $I \subseteq I_0$ there is a string form $J \subseteq I$ such that

$$QT^+ \vdash \forall x \in J\ \forall u,v\ (uEx \ \&\ vEx \rightarrow u=v \lor uEv \lor vEu).$$

(e) For any string form $I \subseteq I_0$ there is a string form $J \subseteq I$ such that

$$QT^+ \vdash \forall y,z \in J\ \forall x\ (xEyz \leftrightarrow xEz \lor x=z \lor \exists w(wEy \ \&\ wz=x)).$$

(f) For any string form $I \subseteq I_0$ there is a string form $J \subseteq I$ such that

$$QT^+ \vdash \forall y,z \in J\ \forall x,x_1,x_2\ (x_1xx_2=yz \rightarrow$$

$$\rightarrow x \subseteq_p y \lor x \subseteq_p z \lor \exists y_1,z_1\ (y_1Ey \ \&\ z_1Bz \ \&\ x=y_1z_1)).$$

(g) For any string concept $I \subseteq I_0$ there is a string concept $J \subseteq I$ such that

$$QT^+ \vdash \forall y,z \in J\ \forall x\ (x \subseteq_p yz \rightarrow x=yz \lor x \subseteq_p y \lor x \subseteq_p z \lor \exists y_1(y_1Ey \ \&\ x=y_1z) \lor$$

$$\lor \exists z_1\ (z_1Bz \ \&\ x=yz_1) \lor \exists y_1,z_1\ (y_1Ey \ \&\ z_1Bz \ \&\ x=y_1z_1)).$$

(h) For any string form $I \subseteq I_0$ there is a string form $J \subseteq I$ such that

$$QT^+ \vdash \forall y,z \in J\ \forall x\ (x \subseteq_p b(yz) \rightarrow$$

$$\rightarrow x=byz \lor x=b \lor x \subseteq_p yz \lor \exists u_2(u_2Byz \ \&\ x=bu_2)).$$

Proof: For (a), see [2], (3.8). For (b) and (c), let $J \equiv I_{LC} \ \&\ I_{LOIS}$. For (d), see [2], (3.10), and then (e) is proved analogously to (b). For (f) take J s in (b), and (g) follows from (b)-(f). Then (h) is obtained as a special case of (g).■

## 4. Coding sequences and pairs of strings by strings

Formalizing recursion requires coding of sequences, and since the kind of recursion used to define the counting functions α and β proceeds on strings, to carry out the formalization of such definitions in concatenation theory we will need to be able to code sequences of strings by strings. The general idea behind the coding goes back to Quine [7], and more recently to Visser [9], but the key for our purposes is to show that the relevant properties of the coding are provable in $QT^+$. We make use of the coding scheme described in [2], pp.86-88 and summarized in [1], §§7-8. (Predicates 'Pref(x,t)', 'Firstf($x,t_1,y,t_2$)', 'Env(t,x)', 'Set(x)' and 'y ε x' are defined there; the formal machinery needed to demonstrate that, modulo the methodology of formula selection, all of the necessary reasoning can indeed be carried out in $QT^+$ is presented in detail in [2], pp.89-263.) In particular, we can establish:

4.1(a) SINGLETON LEMMA. For any string form $I \subseteq I_0$ there is a string form $I_{SNGL} \subseteq I$ such that

$$QT^+ \vdash \forall x \in I_{SNGL}\, \forall u,t_1,t_2\, (Set(x)\ \&\ Firstf(x,t_1,aua,t_2)\ \&\ x=t_1auat_2 \rightarrow$$

$$\rightarrow \forall w\, (w \,\varepsilon\, x \leftrightarrow w=u)).$$

(b) APPENDING LEMMA. For any string form $I \subseteq I_0$ there is a string form $I_{APP} \subseteq I$ such that

$$QT^+ \vdash \forall x,y \in I_{APP}\ \forall t,t_2,t_3(Env(t_2,x)\ \&\ Env(t,y)\ \&\ (t_3a)By\ \&\ Tally_b(t_3)\ \&\ t_2<t_3\ \&$$
$$\&\ \neg\exists u(u\ \varepsilon\ x\ \&\ u\ \varepsilon\ y) \rightarrow \exists z \in I_{APP}\ (Env(t,z)\ \&\ \forall u(u\ \varepsilon\ z \leftrightarrow u\ \varepsilon\ x\ \mathrm{v}\ u\ \varepsilon\ y)).$$

(c) DOUBLETON LEMMA. For any string form $I \subseteq I_0$ there is a string form $I_{DBL} \subseteq I$ such that

$$QT^+ \vdash \forall x \in I_{DBL}\ \forall t_1,t_2,t_3,u,v(Pref(aua,t_1)\ \&\ Pref(ava,t_2)\ \&\ t_1<t_2\ \&\ t_2=t_3\ \&\ u\neq v\ \&$$
$$\&\ x=t_1auat_2avat_3 \rightarrow Set(x)\ \&\ \forall w(w\ \varepsilon\ x \leftrightarrow (w=u\ \mathrm{v}\ w=v)).$$

Proof: See [2], (5.21), (5.46) and (5.58).■

To use the coding of sets to code sequences of strings, we need to populate the coded sets with ordered pairs of arbitrary strings.

Let $\quad Pair[x,y,z] \equiv \exists t \subseteq_p z\ (z=taxatayat\ \&\ MinMax^+T_b(t,xay)).$

(The predicate '$MinMax^+T_b(t,u)$' expressing 't is a shortest non-occurrent b-tally in string u' is defined in [1], §10.) We then have:

4.2 PAIRING LEMMA. (a) For any string form $I \subseteq I_0$ there is a string form $J \subseteq I$ such that

$$QT^+ \vdash \forall x,y \in J\ \exists z \in J\ (Pair[x,y,z]\ \&\ \forall z'(Pair[x,y,z'] \rightarrow z'=z)).$$

(b) For any string form $I \subseteq I_0$ there is a string form $J \subseteq I$ such that

$$QT^+ \vdash \forall z \in J\ \forall x_1,x_2,y_1,y_2\ (Pair[x_1,y_1,z]\ \&\ Pair[x_2,y_2,z] \rightarrow x_1=x_2\ \&\ y_1=y_2).$$

In (a), choose J as in [2], (6.8). For (b), referring to [2], let

$J \equiv I_{3.6}$ & $I_{4.20}$ & $I_{4.23b}$.■

## 5. String recursion

Let p, q be strings, and $f_1$, $f_2$ be functions on strings. Informally, we say that h is defined by string recursion from $f_1$, $f_2$ if

$$h(a) = p \qquad\qquad h(b) = q$$

$$h(y*a) = f_1(y,h(y)) \qquad\qquad h(y*b) = f_2(y,h(y)).$$

We want to justify definitions of this sort in $QT^+$.

Let $I^\lozenge$ be the string form that is the conjunction of the string forms used to obtain the SINGLETON LEMMA, the APPENDING LEMMA, the DOUBLETON LEMMA and the PAIRING LEMMA. The theorem below asserts that, given strings p, q and operations $F_1$, $F_2$ given by formulae satisfying the principal hypothesis, any string form I stronger than $I^\lozenge$ can in turn be strengthened to a string form J containing arbitrarily long length indices for computations of uniquely determined values for successive arguments from J obtained by string recursion from p, q, $F_1$, $F_2$.

STRING RECURSION THEOREM. Let $F_1(y,z,u)$ and $F_2(y,z,u)$ be $\mathcal{L}_C$ formulae, and let $I \subseteq I^\lozenge$ closed under * and downward closed under $\subseteq_p$. Suppose that

$$QT^+ \vdash I(p) \,\&\, I(q),$$

$$QT^+ \vdash \forall y,z \in I\ \exists! u \in I\ F_1(y,z,u), \quad \text{and} \quad QT^+ \vdash \forall y,z \in I\ \exists! u \in I\ F_2(y,z,u).$$

Then there is an $\mathcal{L}_C$ formula $H(y,z)$ and a string form $J \subseteq I$ such that

(i) $QT^+ \vdash \forall y \in J\ \exists! z \in I\ H(y,z)$,

(iia) $QT^+ \vdash \forall y \in I\ (H(a,y) \leftrightarrow y=p)$,

(iib) $QT^+ \vdash \forall y \in I\ (H(b,y) \leftrightarrow y=q)$,

(iiia) $QT^+ \vdash \forall y \in J\ \forall u,z \in I\ (H(y,u) \rightarrow (H(y*a,z) \leftrightarrow F_1(y,u,z)))$,

and (iiib) $QT^+ \vdash \forall y \in J\ \forall u,z \in I\ (H(y,u) \rightarrow (H(y*b,z) \leftrightarrow F_2(y,u,z)))$.

(We read "$\exists! x \in J\ (\ldots)$" as "$\exists x\ (J(x) \,\&\, (\ldots) \,\&\, \forall y (J(y) \,\&\, (\ldots) \rightarrow y=x))$").

Proof: Let Comp(u,m) abbreviate

$Set(u) \,\&\, (a \leq m \rightarrow \exists v \subseteq_p u\ (Pair[a,p,v] \,\&\, v\ \varepsilon\ u)) \,\&$
$\quad \& \, (b \leq m \rightarrow \exists v \subseteq_p u\ (Pair[b,q,v] \,\&\, v\ \varepsilon\ u)) \,\&$

$\& \forall z<m\ \forall u_1,u_2,v_1\ (\text{Pair}[z,u_1,v_1]\ \&\ v_1\ \varepsilon\ u\ \&\ F_1(z,u_1,u_2) \rightarrow$

$\rightarrow \exists v_2 \subseteq_p u\ (\text{Pair}[z*a,u_2,v_2]\ \&\ v_2\ \varepsilon\ u))\ \&$

$\& \forall z<m\ \forall u_1,u_2,v_1\ (\text{Pair}[z,u_1,v_1]\ \&\ v_1\ \varepsilon\ u\ \&\ F_2(z,u_1,u_2) \rightarrow$

$\rightarrow \exists v_2 \subseteq_p u\ (\text{Pair}[z*b,u_2,v_2]\ \&\ v_2\ \varepsilon\ u))\ \&$

$\&\ \forall z,u_1,u_2,v_1,v_2\ (\text{Pair}[z,u_1,v_1]\ \&\ \text{Pair}[z,u_2,v_2]\ \&\ v_1\ \varepsilon\ u\ \&\ v_2\ \varepsilon\ u \rightarrow$

$\rightarrow u_1=u_2\ \&\ v_1=v_2).$

Comp(u,m) means, roughly, that u is a set code for a computation determined by p, q, $F_1,F_2$, in at least m steps where the length indices m are strings ordered by the tree-like ordering $\leq$.

Let MinComp(u,m) abbreviate

$\text{Comp}(u,m)\ \&\ \forall u'\ (\text{Comp}(u',m) \rightarrow \forall y\ (y\ \varepsilon\ u \rightarrow y\ \varepsilon\ u'))\ \&$

$\&\ \forall z,v,w\ (\text{Pair}[z,v,w]\ \&\ w\ \varepsilon\ u \rightarrow (m=a\ \&\ z=a)\ v\ (m=b\ \&\ z=b)\ v$

$v\ \exists n<m\ (z\leq na\ v\ z\leq nb)).$

Let J(m) abbreviate

$I(m)\ \&\ \exists!y\in I\ \exists u\in I\ \exists w\subseteq_p u\ (\text{MinComp}(u,m)\ \&\ \text{Pair}[m,y,w]\ \&\ w\ \varepsilon\ u).$

Finally, let H(m,y) abbreviate

$\exists u,w\ (\text{MinComp}(u,m)\ \&\ \text{Pair}[m,y,w]\ \&\ w\ \varepsilon\ u).$

For detailed verification that J and H have the desired properties see Appendix.■

We are now in the position to define the counting functions $\alpha$ and $\beta$.

Let $p = bb$, $q = b$, $F_1(y,z,u) \equiv y=y \,\&\, z=Su$ and $F_2(y,u,z) \equiv y=y \,\&\, z=u$.

Then the principal hypothesis of the String Recursion Theorem holds trivially.

Applying the Theorem we obtain a formula $A^{\#}(y,z)$ and a string form $I_\alpha \subseteq I$ such that

$(i^\alpha)$ $\quad QT^+ \vdash \forall y \in I_\alpha \, \exists! z \in I \;\; A^{\#}(y,z)$,

$(iia^\alpha)$ $\quad QT^+ \vdash \forall z \in I \, (A^{\#}(a,z) \leftrightarrow z=bb)$,

$(iib^\alpha)$ $\quad QT^+ \vdash \forall z \in I \, (A^{\#}(b,z) \leftrightarrow z=b)$,

$(iiia^\alpha)$ $\quad QT^+ \vdash \forall y \in I_\alpha \, \forall u,z \in I \, (A^{\#}(y,u) \;\rightarrow\; (A^{\#}(y*a,z) \rightarrow z=u*b))$,

$(iiib^\alpha)$ $\quad QT^+ \vdash \forall y \in I_\alpha \, \forall u,z \in I \, (A^{\#}(y,u) \;\rightarrow\; (A^{\#}(y*b,z) \rightarrow z=u))$.

Informally speaking, $A^{\#}(y,z)$ defines the graph of the function $\alpha$.

Exactly analogously, by letting p and q, and $F_1$, $F_2$, respectively, exchange places, we apply the Theorem to obtain a formula $B^{\#}(y,z)$ defining the graph of the function $\beta$ and a string form $I_\beta \subseteq I$ such that

$(i^\beta)$ $\quad QT^+ \vdash \forall y \in I_\beta \, \exists! z \in I \;\; B^{\#}(y,z)$,

$(iia^\beta)$ $\quad QT^+ \vdash \forall z \in I \, (B^{\#}(a,z) \leftrightarrow z=b)$,

$(iib^\beta)$ $\quad QT^+ \vdash \forall z \in I \, (B^{\#}(b,z) \leftrightarrow z=bb)$,

$(iiia^\beta)$ $\quad QT^+ \vdash \forall y \in I_\beta \, \forall u,z \in I \, (B^{\#}(y,u) \;\rightarrow\; (B^{\#}(y*a,z) \rightarrow z=u))$,

$(iiib^\beta)$ $\quad QT^+ \vdash \forall y \in I_\beta \, \forall u,z \in I \, (B^{\#}(y,u) \;\rightarrow\; (B^{\#}(y*b,z) \rightarrow z=u*b))$.

We can then prove that $\alpha$ and $\beta$ correctly count b's in b-tallies:

5.1(a) For any string form $I \subseteq I_0$ there is a string form $J \subseteq I$ such that

$$QT^+ \vdash A^\#(a,Sb) \,\&\, \forall x \in J\ \forall y \in I\ (Tally_b(x) \,\&\, A^\#(x,y) \rightarrow y=b) \quad \text{and}$$

$$QT^+ \vdash \forall x \in J\ \forall y \in I\ (Tally_a(x) \,\&\, B^\#(x,y) \rightarrow y=b).$$

(I.e., '$\alpha(a)=1$' and '$Tally_b(x) \rightarrow \alpha(x)=0$', and '$Tally_a(x) \rightarrow \beta(x)=0$'.)

(b) For any string form $I \subseteq I_0$ there is a string form $J \subseteq I$ such that

$$QT^+ \vdash \forall x \in J\ \forall y \in I\ (Tally_b(x) \,\&\, B^\#(x,y) \rightarrow y=x*b).$$

Informally, $Tally_b(x) \rightarrow \beta(x)=length(x)$.

We now verify that the functions $\alpha$ and $\beta$ are indeed additive. Let $I_{Add}$ be as in 3.5(a).

5.2(a) For any string form $I \subseteq I_\alpha$ and $I \subseteq I_{Add}$ there is a string form $J \equiv I_{Add\alpha} \subseteq I$ such that

$$QT^+ \vdash \forall x,y \in J\ \forall u,v,w\ (A^\#(x,u) \,\&\, A^\#(y,v) \,\&\, AddTally(u,v,w) \rightarrow A^\#(x*y,w)).$$

("$\alpha(x*y) = \alpha(x) + \alpha(y)$")

(b) For any string form $I \subseteq I_\beta$ and $I \subseteq I_{Add}$ there is a string form $J \equiv I_{Add\beta} \subseteq I$ such that

$$QT^+ \vdash \forall x,y \in J\ \forall u,v,w\ (B^\#(x,u) \,\&\, B^\#(y,v) \,\&\, AddTally(u,v,w) \rightarrow B^\#(x*y,w)).$$

("$\beta(x*y) = \beta(x) + \beta(y)$")

<u>Proof</u>: See Appendix.■

## 6. Formal Construction of the Interpretation

Let $I_{Add\alpha}$ be the string form obtained from $I_0$ by the series of modifications described in §§3-5 up to and including 5.2(a). Analogously for and $I_{Add\beta}$ and 5.2(b).

Let $J^* \equiv I_{Add\alpha} \,\&\, I_{Add\beta}$.

Then $J^* \subseteq I_{Add\alpha}$ and $J^* \subseteq I_{Add\beta}$ and $J^* \subseteq I_{Add}$ as well as $J^* \subseteq I^{\lozenge}$. We may also assume that $J^*$ is closed under *, and downward closed under $\leq$ and $\subseteq_p$. Hence it may be assumed that the string form $J^*$ is also closed under Addtally and the functions $\alpha$ and $\beta$.

We then formally define Æ(x) as

$\exists y,z\ (A^\#(x,y) \,\&\, B^\#(x,z) \ \&\ y=Sz)\ \&$

$\&\ \forall u,v,w\ (uBx \,\&\, A^\#(u,v) \,\&\, B^\#(u,w) \ \rightarrow\ v\leq w).$

(These are conditions (c1)-(c2) from §1.)

We set $I^*(x) \ \equiv$ Æ(x) & $J^*(x)$.

The formula $I^*(x)$ will formally define in $QT^+$ the domain of interpretation of theory T. We now proceed to formally verify the translations of the axioms of T by derivations in $QT^+$.

6.1(a) $QT^+ \vdash I^*(x) \,\&\, x_2Ex \ \rightarrow\ \forall u,v\ (A^\#(x_2,u) \,\&\, B^\#(x_2,v) \ \rightarrow\ Sv\leq u).$

(b) $QT^+ \vdash I^*(x) \,\&\, I^*(y) \,\&\, z=bxy \ \rightarrow\ I^*(z).$

(c) $QT^+ \vdash I^*(x) \,\&\, I^*(u) \,\&\, bxy=buv \;\rightarrow\; x=u \,\&\, y=v.$

(d) $QT^+ \vdash I^*(x) \;\rightarrow (x \subseteq_p a \leftrightarrow x=a).$

(e) $QT^+ \vdash I^*(x) \,\&\, I^*(y) \,\&\, I^*(z) \;\rightarrow\; (x \subseteq_p byz \leftrightarrow x=byz \;v\; x \subseteq_p y \;v\; x \subseteq_p z).$

Proof: See Appendix. We give the details of the proof of (e) to illustrate the flavor of the type of formal argument used.

Assume $M \vDash x \subseteq_p byz$ where $M \vDash I^*(x) \,\&\, I^*(y) \,\&\, I^*(z)$.

Then $M \vDash J^*(x) \,\&\, J^*(y) \,\&\, J^*(z)$ and also $M \vDash Æ(x) \,\&\, Æ(y) \,\&\, Æ(z)$.

By $(i^\alpha)$ and $(i^\beta)$, $M \vDash \exists! x_1 \in J^* \, A^\#(x,x_1) \,\&\, \exists! x_2 \in J^* \, B^\#(x,x_2)$.

From $M \vDash x \subseteq_p byz$ by 3.7(h) we have that

$$M \vDash x=byz \;v\; x=b \;v\; x \subseteq_p yz \;v\; \exists u(uByz \,\&\, x=bu).$$

We distinguish the cases:

(1) $M \vDash \exists u(uByz \,\&\, x=bu)$.

Then by (QT2), $M \vDash x \neq a$. From $M \vDash I^*(y) \,\&\, I^*(z)$, by 6.1(b), $M \vDash I^*(byz)$.

From $M \vDash uByz$, $M \vDash \exists v \; uv=yz$, hence $M \vDash b(uv)=b(yz)$, also $M \vDash (bu)v=b(yz)$. Thus $M \vDash buBb(yz)$, hence $M \vDash xBb(yz)$.

From $M \vDash I^*(byz)$, $M \vDash Æ(byz)$, whence $M \vDash x_1 \leq x_2$. But from $M \vDash Æ(x)$, $M \vDash x_1=Sx_2$, and we have $M \vDash x_1 \leq x_2 < Sx_2 = x_1$, contradicting $M \vDash I_0(x_1)$.

Hence (1) is ruled out.

(2) $M \vDash x=b$.

Then by (QT2), $M \vDash x \neq a$, and from $M \vDash Æ(x)$, we have $M \vDash bBx$. But then $M \vDash bBb$, contradicting (QT2). Hence (2) is also ruled out.

(3) $M \vDash x \subseteq_p yz$.

By 3.7(g), $M \models x{=}yz \ v \ x \subseteq_p y \ v \ x \subseteq_p z \ v \ \exists y_1(y_1 E y \ \& \ x{=}y_1 z) \ v$

$v \ \exists z_1(z_1 B z \ \& \ x{=}yz_1) \ v \ \exists y_1,z_1(y_1 E y \ \& \ z_1 B z \ \& \ x{=}y_1 z_1)$.

(3i) $M \models x{=}yz$.

By ($i^\alpha$) and ($i^\beta$), $M \models \exists! y_1 \in J^* \ A^\#(y,y_1) \ \& \ \exists! y_2 \in J^* \ B^\#(y,y_2)$,

and further $M \models \exists! z_1 \in J^* \ A^\#(z,z_1) \ \& \ \exists! z_2 \in J^* \ B^\#(z,z_2)$.

From $M \models Æ(y)$, $M \models y_1{=}Sy_2$, and from $M \models Æ(z)$, $M \models z_1{=}Sz_2$.

By 3.5(a), $M \models \exists! p_1 \in J^*(Tally_b(p_1) \ \& \ Addtally(y_1,z_1,p_1))$

and $M \models \exists! p_2 \in J^*(Tally_b(p_2) \ \& \ Addtally(y_2,z_2,p_2))$.

Then from $M \models A^\#(y,y_1) \ \& \ A^\#(z,z_1)$, by 5.2(a), $M \models A^\#(y{*}z,p_1)$, and

from $M \models B^\#(y,y_2) \ \& \ B^\#(z,z_2)$, by 5.2(b), $M \models B^\#(y{*}z,p_2)$,

$\Rightarrow$ from $M \models y_1{=}Sy_2 \ \& \ z_1{=}Sz_2$, $M \models Addtally(Sy_2,Sz_2,p_1)$.

On the other hand, from $M \models Addtally(y_2,z_2,p_2)$, by 3.4(e),

$M \models Addtally(y_2,Sz_2,Sp_2)$, whence by 3.5(d), $M \models Addtally(Sy_2,Sz_2,SSp_2)$.

By single-valuedness of Addtally, we then have $M \models p_1{=}SSp_2$.

From hypothesis $M \models x{=}yz \ \& \ A^\#(y{*}z,p_1) \ \& \ B^\#(y{*}z,p_2)$,

$M \models A^\#(x,p_1) \ \& \ B^\#(x,p_2)$.

Hence from $M \models A^\#(x,x_1) \ \& \ B^\#(x,x_2)$, by single-valuedness of $A^\#$ and $B^\#$,

$M \models p_1{=}x_1 \ \& \ p_2{=}x_2$.

Thus from $M \models p_1{=}SSp_2$, we have $M \models x_1{=}SSx_2$. But from $M \models Æ(x)$ we have $M \models x_1{=}Sx_2$, whence $M \models x_1{=}Sx_1$. But then from $M \models x_1{<}Sx_1$, we obtain $M \models x_1{<}x_1$, contradicting $M \models I_0(x_1)$. Hence (3i) is ruled out.

(3ii) $M \models \exists y_1(y_1 E y \ \& \ x{=}y_1 z)$.

Then $M \models y_1 Bx$.

By ($i^\alpha$) and ($i^\beta$), $M \models \exists! u_1 \in J^* A^\#(y_1,u_1)$ & $\exists! u_2 \in J^* B^\#(y_2,u_2)$.

From $M \models Æ(x)$ & $y_1Bx$, $M \models u_1 \leq u_2$, whereas from $M \models I^*(y)$ & $y_1Ey$, by 6.1(a), $M \models Su_2 \leq u_1$. But then $M \models u_2 < Su_2 \leq u_2$, contradicting $M \models I_0(u_2)$. This rules out (3ii).

(3iii) $M \models \exists z_1(z_1Bz$ & $x=yz_1)$.

Then $M \models yBx$. By ($i^\alpha$) and ($i^\beta$), $M \models \exists! y_1 \in J^* A^\#(y,y_1)$ & $\exists! y_2 \in J^* B^\#(y,y_2)$. From $M \models Æ(x)$ & $yBx$, $M \models y_1 \leq y_2$. But from $M \models Æ(y)$, $M \models y_1 = Sy_2$, and we obtain $M \models y_1 \leq y_2 < Sy_2 = y_1$, contradicting $M \models I_0(y_1)$. Hence (3iii) is ruled out.

(3iv) $M \models \exists y_1,z_1(y_1Ey$ & $z_1Bz$ & $x=y_1z_1)$.

This is ruled out by reasoning as in either (3ii) or (3iii).

We then conclude under the principal hypothesis that

$$M \models x \subseteq_p yz \rightarrow x \subseteq_p y \ \text{v} \ x \subseteq_p z$$

and further that $M \models x \subseteq_p byz \rightarrow x = byz \ \text{v} \ x \subseteq_p y \ \text{v} \ x \subseteq_p z$.

The converse is immediate from the definition of $\subseteq_p z$.■

Taking the formula Æ(x) from §6 to define the domain, and interpreting the non-logical vocabulary $\mathcal{L}_T = \{0, (\ ), \sqsubseteq\}$ of T by *a*, *b*xy and $\subseteq_p$, resp., as explained in §2, we have that 6.1(b)-(e), along with the fact that $QT^+ \vdash bxy \neq a$, suffice to establish formal interpretability of T in $QT^+$. On the other hand, from [1], building on previous work of Halpern and Collins, Wilkie, Visser, Grzegorczyk and Ganea, we have that

$$\mathrm{TC} \equiv_I \mathrm{QT}^+ \equiv_I \mathrm{AST} \equiv_I \mathrm{AST{+}EXT} \equiv_I \mathrm{Q} \equiv_I$$

Since, by [6], $\mathrm{Q} \leq_I \mathrm{T}$, this suffices to establish

WEAK ESSENTIALLY UNDECIDABLE THEORIES: FIRST MUTUAL INTERPRETABILITY THEOREM.

$$\mathrm{T} \equiv_I \mathrm{QT}^+ \equiv_I \mathrm{QT}_0 \equiv_I \mathrm{TC} \equiv_I \mathrm{Q} \equiv_I \mathrm{AST}.$$

In addition, each of the theories above is mutually interpretable with AST+EXT, and Buss's theory $S^1_2$ (see Ferreira and Ferreira [4]).

## §7. R and its variants

We now consider the expanded vocabulary $\mathcal{L}_{C,\sqsubseteq^*} = \{a, b, *, \sqsubseteq^*\}$ with two individual constants – the digits *a*, *b* – a single binary operation symbol * and a 2-place relational symbol $\sqsubseteq^*$. Each variable-free term of $\mathcal{L}_{C,\sqsubseteq^*}$ represents a finite string of *a*'s and/or *b*'s, and each such string may have multiple variable-free terms as its representations, differing in the arrangement of parentheses indicating the order of applications of the term operation *.

Recalling the theory WT described in the introduction formulated in $\mathcal{L}_T = \{0, (\,), \sqsubseteq\}$, we are going to single out $\mathcal{L}_{C,\sqsubseteq^*}$ terms that represent tree-like strings obtained from variable-free terms of $\mathcal{L}_T$ as described in §1. With each variable-free term v of $\mathcal{L}_T$ we associate a unique $\mathcal{L}_{C,\sqsubseteq^*}$ term $v^\tau$ as follows:

$$0^\tau \equiv a \qquad\qquad (u,v)^\tau \equiv b*(u^\tau * v^\tau).$$

The $\mathcal{L}_{QT,\sqsubseteq^*}$ term $v^\tau$ is an Æ string that codes v.

If $\mathcal{S}(v)$ is the set of all (variable-free) subterms of v, let

$$\Sigma(t) = \{\, u^\tau \mid \text{for some } \mathcal{L}_T\text{-term } v,\ \ u \in \mathcal{S}(v) \text{ and } t{=}v^\tau\}.$$

We then let $\Sigma^\tau = \bigcup_{v \in \mathcal{S}} \Sigma\,(v^\tau)$, where $\mathcal{S}$ is the set all variable-free terms of $\mathcal{L}_T$. A straightforward induction on the complexity of $\mathcal{L}_T$ terms establishes that the mapping $\tau$ is 1-1.

Let WQT be the first-order theory formulated in $\mathcal{L}_{C,\sqsubseteq^*}$ with the following axioms:

(WQT1) $\neg(s{=}t)$ for any distinct terms $s, t \in \Sigma^\tau$,

(WQT2) $\forall z\, (z \sqsubseteq^* b*(s*t) \leftrightarrow z{=}b*(s*t) \ \text{v}\ z \sqsubseteq^* s \ \text{v}\ z \sqsubseteq^* t)$

for all terms $s, t \in \Sigma^\tau$,

(WQT3) $\forall z\,(z \sqsubseteq^* a \leftrightarrow z=a)$.

Here, (WQT1) and (WQT2) are axiom schemas with infinitely many instances.

We now define a formal interpretation $(\tau)$ of WT in WQT. Let the formula

$$T^*(x) \equiv x=a \vee \exists y,z\; x=b*(y*z)$$

define the domain. Interpret 0 by $a$, the binary term building operation ( , ) of $\mathcal{L}_T$ by $b*(x*y)$, and $\sqsubseteq$ by $\sqsubseteq^*$. We then have immediately:

$$\text{WQT} \vdash T^*(0^\tau),$$

$$\text{WQT} \vdash T^*(y) \;\&\; T^*(z) \rightarrow T^*(b*(y*z)).$$

A trivial induction on the complexity of $\mathcal{L}_T$ terms verifies that each $v \in \mathcal{S}$ is interpreted by $v^\tau \in \Sigma^\tau$ in WQT. Since the map $\tau$ is 1-1, we have that

$$\text{WQT} \vdash (\neg(u=v))^{(\tau)},$$

$\neg(u^\tau = v^\tau)$ being the translations $(\neg(u=v))^{(\tau)}$ of the instances of axiom schema (WT1) of WT, for distinct $u, v \in \mathcal{S}$.

Consider now an instance of the schema (WT2), for some $v \in \mathcal{S}$:

$$\forall x\, (x \sqsubseteq v \leftrightarrow \bigvee_{u \in \mathcal{S}(v)} x=u).$$

If v is the atomic term 0, we have that $\mathcal{S}(v) = \{0\}$. Hence the formula in question is $\forall x\, (x \sqsubseteq 0 \leftrightarrow x=0)$.

But, by (WQT3), $\mathrm{WQT} \vdash \forall x\, (x \sqsubseteq^{*} a \leftrightarrow x=a)$.

Hence, *a fortiori*, $\mathrm{WQT} \vdash \forall x\, (T^{*}(x) \rightarrow (x \sqsubseteq^{*} a \leftrightarrow x=a))$, which is the

(τ)-translation of the above instance of (WT2).

Consider now $t \in \mathcal{S}$ of the form (u,v). Note that

$$\mathcal{S}(t) = \mathcal{S}(u) \cup \mathcal{S}(v) \cup \{t\}.$$

Hence $\Sigma(t^{\tau}) = \Sigma(u^{\tau}) \cup \Sigma(v^{\tau}) \cup \{t^{\tau}\}$. (†).

Assume now that

$\mathrm{WQT} \vdash [\forall z\, (z \sqsubseteq u \leftrightarrow \bigvee_{s \in \mathcal{S}(u)} z=s)]^{(\tau)}$ and $\mathrm{WQT} \vdash [\forall z\, (z \sqsubseteq v \leftrightarrow \bigvee_{s \in \mathcal{S}(v)} z=s)]^{(\tau)}$.

Then $\mathrm{WQT} \vdash \forall z\, (T^{*}(z) \rightarrow (z \sqsubseteq^{*} u^{\tau} \leftrightarrow \bigvee_{s^{\tau} \in \Sigma(u^{\tau})} z=s^{\tau}))$

and $\mathrm{WQT} \vdash \forall z\, (T^{*}(z) \rightarrow (z \sqsubseteq^{*} v^{\tau} \leftrightarrow \bigvee_{s^{\tau} \in \Sigma(v^{\tau})} z=s^{\tau}))$.

Let M be a model of WQT. Assume $M \vDash T^{*}(x)$ and consider $M \vDash x \sqsubseteq^{*} t$.

We have that $t^{\tau}$ is in fact $b^{*}(u^{\tau} * v^{\tau})$. Hence

$\Leftrightarrow$ $M\vDash x\sqsubseteq^{*}b*(u^{\tau}*v^{\tau})$ $\Leftrightarrow$ by (WQT2), $M\vDash x=b*(u^{\tau}*v^{\tau})$ v $x\sqsubseteq^{*}u^{\tau}$ v $x\sqsubseteq^{*}v^{\tau}$ $\Leftrightarrow$

$\Leftrightarrow$ $M\vDash x=b*(u^{\tau}*v^{\tau})$ v $\bigvee_{s^{\tau}\in\Sigma(u^{\tau})} x=s^{\tau}$ v $\bigvee_{s^{\tau}\in\Sigma(v^{\tau})} x=s^{\tau}$ $\Leftrightarrow$

$\Leftrightarrow$ $M\vDash \bigvee_{s^{\tau}\in\Sigma(t^{\tau})} x=s^{\tau}$

using (†). Therefore,

$$\text{WQT}\vdash\forall x\,(T^{*}(x)\to(x\sqsubseteq^{*}t^{\tau}\leftrightarrow\bigvee_{s^{\tau}\in\Sigma(t^{\tau})}x=s^{\tau})),$$

that is, $\quad\text{WQT}\vdash[\forall x\,(x\sqsubseteq t\leftrightarrow\bigvee_{s\in\mathcal{S}(t)}x=s)]^{(\tau)}$.

Hence the (τ)-translation of each instance of (WT2) is also provable in WQT.

We conclude that

7.1. WT $\leq_{I}$ WQT.

The theory WQT is not recognizably a concatenation theory: the axioms make no substantive assumptions about the binary operation *, not even associativity. On that account, it might be considered at best as a "pseudo-concatenation" notational variant of WT. We now consider another first-order theory, WQT*, formulated in the same vocabulary $\mathcal{L}_{C,\sqsubseteq^{*}}=\{a, b, *, \sqsubseteq^{*}\}$ as WQT, with the following axioms: for each variable-free term t of $\mathcal{L}_{C,\sqsubseteq^{*}}$,

(WQT*1) $\forall x,y,z\,(x*(y*z)\sqsubseteq_{p}t$ v $(x*y)*z\sqsubseteq_{p}t\to x*(y*z)=(x*y)*z)$

(WQT*2) $\forall x,y\ (x{*}y \sqsubseteq_p t \ \rightarrow \neg(x{*}y=a)\ \&\ \neg(x{*}y=b))$

(WQT*3) $\forall x,y\ ((a{*}x \sqsubseteq_p t\ \&\ a{*}y \sqsubseteq_p t \rightarrow (a{*}x=a{*}y \rightarrow x=y))\ \&$

$\&\ (b{*}x \sqsubseteq_p t\ \&\ b{*}y \sqsubseteq_p t \rightarrow (b{*}x=b{*}y \rightarrow x=y))\ \&$

$\&\ (x{*}a \sqsubseteq_p t\ \&\ y{*}a \sqsubseteq_p t \rightarrow (x{*}a=y{*}a \rightarrow x=y))$

$\&\ (x{*}b \sqsubseteq_p t\ \&\ x{*}y \sqsubseteq_p t \rightarrow (x{*}b=y{*}b \rightarrow x=y)))$

(WQT*4) $\forall x,y\ ((a{*}x \sqsubseteq_p t\ \&\ b{*}y \sqsubseteq_p t \ \rightarrow\ \neg(a{*}x=b{*}y))\ \&$

$\&\ (x{*}a \sqsubseteq_p t\ \&\ y{*}b \sqsubseteq_p t \ \rightarrow\ \neg(x{*}a=y{*}b)))$

(WQT*5) $\forall x \sqsubseteq_p t\ (x=a\ \ v\ \ x=b\ \ v\ \ ((aBx\ v\ bBx)\ \&\ (aEx\ v\ bEx)))$

(WQT*6) $\forall y,z\ (b{*}(y{*}z) \sqsubseteq^{*} t \ \rightarrow$

$\rightarrow\ \forall x\ (x \sqsubseteq^{*} b{*}(y{*}z) \ \leftrightarrow\ x=b{*}(y{*}z)\ \ v\ \ x \sqsubseteq^{*} y\ \ v\ \ x \sqsubseteq^{*} z))$

(WQT*7) $\forall z\ (z \sqsubseteq^{*} a \ \leftrightarrow z=a)$

(WQT*8) $\forall x,y\ (x \sqsubseteq^{*} y\ \&\ y \sqsubseteq^{*} x \ \rightarrow\ x=y)$

(WQT*9) $\forall x,y\ (x \sqsubseteq^{*} y\ \&\ y \sqsubseteq^{*} z \rightarrow\ y \sqsubseteq^{*} z)$

Here we use the following abbreviations:

$xBy \equiv \exists z\ y=x{*}z \qquad\qquad xEy \equiv \exists z\ y=z{*}x,$

and $x \sqsubseteq_p y \ \equiv\ x=y \ \text{v}\ xBy \ \text{v}\ xEy \ \text{v}\ \exists z_1,z_2\ y=z_1*(x*z_2)\ \text{v}\ \exists z_1,z_2\ y=(z_1*x)z_2$.

Then $\forall x \sqsubseteq_p u\ \varphi \ \equiv\ \forall x\,(x \sqsubseteq_p u \ \rightarrow\ \varphi)$, where x does not occur in the term u.

Also, $\forall x \sqsubseteq^* u\ \varphi \ \equiv\ \forall x\,(x \sqsubseteq^* u \ \rightarrow\ \varphi)$.

($WQT^*1$)-($WQT^*6$) are axiom schemas with infinitely many instances, one for each variable-free term t. The schemas ($WQT^*1$)-($WQT^*5$) are “bounded” versions of the axioms (QT1)-(QT5) of $QT^+$. Schema ($WQT^*6$) is a “bounded” generalization of schema (WQT2) of WQT. In light of that, $WQT^*$ may be naturally interpreted as a hybrid basic theory of finite strings and trees: the intended domain are the finite strings of *a*’s and/or *b*’s, * is interpreted as the concatenation operation, and $\sqsubseteq^*$ as the substring relation between Æ strings.

$WQT^*$ is an extension of WQT. First, we note the following:

7.2. For any distinct terms $s, t \in \Sigma^\tau$, $WQT^* \vdash \neg(s=t)$.

<u>Proof</u>: We argue by (meta-theoretic) induction on the number of digits in s, t. If either one of s or t is the single digit a, this is immediate by ($WQT^*2$). If neither s nor t are single digits, let $s_1 \ldots s_m$ and $t_1 \ldots t_n$ be their successive digits (ignoring parentheses), and let $s_i \neq t_i$ be the leftmost digit where s and t differ. Then $s = s_1 \ldots s_{i-1} s_i s^-$ and $t = t_1 \ldots t_{i-1} t_i t^-$ where $s^- = s_{i+1} \ldots s_m$ and $t^- = t_{i+1} \ldots t_n$. By (WQT*4), $WQT^* \vdash \neg(s_i s^- = t_i t^-)$. By repeatedly applying (WQT*1) and (WQT*3) we obtain $WQT^* \vdash \neg(s_1 \ldots s_{i-1} s_i s^- = t_1 \ldots t_{i-1} t_i t^-)$, that is, $WQT^* \vdash \neg(s=t)$, as required.■

Hence in particular all instances of schema (WQT1) are provable in WQT*. Consider an instance of (WQT2) for terms $s, t \in \Sigma^\tau$,

$$\forall z\, (z \sqsubseteq^* b^*(s^*t) \leftrightarrow z = b^*(s^*t) \lor z \sqsubseteq^* s \lor z \sqsubseteq^* t)$$

Now, $WQT^* \vdash b^*(s^*t) = b^*(s^*t)$, so $WQT^* \vdash b^*(s^*t) \sqsubseteq_p b^*(s^*t)$. From (WQT*6),

$$WQT^* \vdash \forall x\, (x \sqsubseteq^* b^*(s^*t) \leftrightarrow x = b^*(s^*t) \lor x \sqsubseteq^* s \lor x \sqsubseteq^* t).$$

Hence each instance of (WQT2) is provable in WQT*. Given that (WQT3) is (WQT*7), this is enough to establish that WQT* is an extension of WQT.

On the other hand, we also have:

7.3. WQT* is locally finitely satisfiable.

That is, each finite subset of its non-logical axioms has a finite model.

Proof: See Appendix.■

By Visser's Theorem, it follows that $WQT^*$ is interpretable in R.

Since by [6], $R \leq_I WT$, we then have

WEAK ESSENTIALLY UNDECIDABLE THEORIES: SECOND MUTUAL INTERPRETABILITY THEOREM. $R \equiv_I WTC^{-\varepsilon} \equiv_I WT \equiv_I WQT \equiv_I WQT^*$.

For definition of the theory $WTC^{-\varepsilon}$ see [5].

R E F E R E N C E S

# A P P E N D I X

2.2 $\Sigma^* \vDash Æ(x) \leftrightarrow x=a \vee \exists! y,z\,(x=b(yz)\ \&\ Æ(y)\ \&\ Æ(z))$.

Proof: (⇐) Assume $\Sigma^* \vDash Æ(y)\ \&\ Æ(z)\ \&\ x=byz$. Then

$$\Sigma^* \vDash \alpha(y) = \beta(y)+1\ \&\ \alpha(z) = \beta(z)+1.$$

Now, $\Sigma^* \vDash \alpha(x) = \alpha(byz) = \alpha(yz) = \alpha(y)+\alpha(z)$

and $\Sigma^* \vDash \beta(x) = \beta(byz) = \beta(b)+\beta(yz) = \beta(y)+\beta(z)+1$. Then

$\Sigma^* \vDash \alpha(x) = \alpha(y)+\alpha(z) = (\beta(y)+1)+(\beta(z)+1) = (\beta(y)+\beta(z)+1)+1 = \beta(x)+1$

which verifies (c1). For (c2), assume $\Sigma^* \vDash uBx$, i.e., $\Sigma^* \vDash uBbyz$.

Then $\Sigma^* \vDash u=b \vee uBby \vee u=by \vee \exists z_1(z_1Bz\ \&\ u=byz_1)$.

If (a) $\Sigma^* \vDash u=b$, then $\Sigma^* \vDash \alpha(u) = \alpha(b) = 0 < 1 = \beta(b) = \beta(u)$.

If (b) $\Sigma^* \vDash uBby$, then $\Sigma^* \vDash \exists y_1(y_1By\ \&\ u=by_1)$. Then from $\Sigma^* \vDash Æ(y)$,

$\Sigma^* \vDash \alpha(y_1) \leq \beta(y_1)$, whence $\Sigma^* \vDash \alpha(u) = \alpha(by_1) = \alpha(y_1)$ and

$\Sigma^* \vDash \beta(u) = \beta(by_1) = \beta(b)+\beta(y_1) = \beta(y_1)+1$.

Hence $\Sigma^* \vDash \alpha(u) = \alpha(y_1) \leq \beta(y_1) < \beta(y_1)+1 = \beta(u)$.

Suppose (c) $\Sigma^* \vDash u{=}by$. Then from $\Sigma^* \vDash Æ(y)$, $\Sigma^* \vDash \alpha(y) = \beta(y)+1$, and we have $\Sigma^* \vDash \alpha(u) = \alpha(by) = \alpha(y)$ and

$\Sigma^* \vDash \beta(u) = \beta(by) = \beta(b)+\beta(y) = \beta(y_1)+1.$

Hence $\Sigma^* \vDash \alpha(u) = \alpha(y) = \beta(y)+1 = \beta(u)$, so $\Sigma^* \vDash \alpha(u) \leq \beta(u)$.

Finally, suppose (d) $\Sigma^* \vDash \exists z_1(z_1By \ \& \ u{=}byz_1)$. Then from $\Sigma^* \vDash Æ(z)$,

$\Sigma^* \vDash \alpha(z_1) \leq \beta(z_1)$, and from $\Sigma^* \vDash Æ(y)$, $\Sigma^* \vDash \alpha(y) = \beta(y)+1$. Then

$\Sigma^* \vDash \alpha(u) = \alpha(byz_1) = \alpha(yz_1) = \alpha(y)+\alpha(z_1)$ and

$\Sigma^* \vDash \beta(u) = \beta(byz_1) = \beta(b)+\beta(yz_1) = \beta(y)+\beta(z_1)+1.$ Hence

$$\Sigma^* \vDash \alpha(u) = \alpha(y)+\alpha(z_1) = (\beta(y)+1)+\alpha(z_1) \leq$$
$$\leq (\beta(y)+1)+\beta(z_1) = \beta(y)+\beta(z_1)+1 = \beta(u).$$

Thus $\Sigma^* \vDash \alpha(u) \leq \beta(u)$. This completes the proof of (c2). So $\Sigma^* \vDash Æ(x)$.■

3.5(i) For any string form $I \subseteq I_0$ there is a string form $J \subseteq I$ such that

$$QT^+ \vdash \forall x_2,y_1,y_2 \in J\ \forall x_1,z_1,z_2\ (Tally_b(x_2)\ \&\ Tally_b(y_1)\ \&\ Tally_b(y_2)\ \&$$
$$\&\ Addtally(x_1,x_2,z_1)\ \&\ Addtally(y_1,y_2,z_2)\ \&\ x_1 \leq y_1\ \&\ z_1 = Sz_2 \rightarrow$$
$$\rightarrow Sy_2 \leq x_2).$$

Let $J(y) \equiv I_{LC}$ & $I_{CTC}$ & $I_{6.2(a)}$ & $I_{COMM}$.

Assume $M \vDash Addtally(x_1,x_2,z_1)$ & $Addtally(y_1,y_2,z_2)$

where $M \vDash x_1 \leq y_1$ & $z_1=Sz_2$ and $M \vDash Tally_b(x_2)$ & $Tally_b(y_1)$ & $Tally_b(y_2)$ and $M \vDash J(y_1)$.

From $M \vDash Tally_b(y_1)$ & $x_1 \leq y_1$, $M \vDash Tally_b(x_1)$.

By (f), $M \vDash \exists u_1(Tally_b(u_1)$ & $Addtally(u_1,x_1,y_1))$, whereas by 3.5(a),

$M \vDash \exists! p_1 \in J\ Addtally(u_1,x_1,p_1)$.

Then, by single-valuedness of Addtally, $M \vDash y_1=p_1$, whence from hypothesis $M \vDash Addtally(y_1,y_2,z_2)$, $M \vDash Addtally(p_1,y_2,z_2)$.

On the other hand, $M \vDash \exists! p_2 \in J\ Addtally(x_1,u_1,p_2)$, whence by (g),

$M \vDash Addtally(u_1,x_1,p_2)$.

But then from $M \vDash Addtally(u_1,x_1,p_1)$, by single-valuedness of Addtally,

$$M \vDash p_1=p_2.$$

Hence $M \vDash Addtally(p_2,y_2,z_2)$, and from 3.4(e) we obtain

$$M \vDash Addtally(p_2,y_2{*}b,z_2{*}b).$$

From $M \vDash Addtally(y_1,y_2,z_2)$, by 3.5(a), $M \vDash Tally_b(z_2)$.

So from $M \vDash Tally_b(y_2)$ & $Tally_b(z_2)$, $M \vDash Addtally(p_2,Sy_2,Sz_2)$.

Again by 3.5(a),

$$M \vDash \exists! v_1 \in J\ Addtally(u_1,Sy_2,v_1) \text{ and } M \vDash \exists! w_1 \in J\ Addtally(x_1,v_1,w_1).$$

From $M \vDash Addtally(x_1,u_1,p_2)$ by (h), $M \vDash Sz_2=w_1$.

From hypothesis $M \vDash z_1=Sz_2$, $M \vDash z_1=w_1$, so $M \vDash Addtally(x_1,v_1,z_1)$.

On the other hand, from hypothesis $M \vDash Addtally(x_1,x_2,z_1)$, by (e),

$$M \vDash x_2 = v_1,$$

Hence from $M \vDash Addtally(u_1,Sy_2,v_1)$, $M \vDash Addtally(u_1,Sy_2,x_2)$.

But then from $M \vDash Tally_b(u_1)$, by (h), $M \vDash Sy_2 \leq x_2$, as required.■

Let $$MaxT_b(t,w) \equiv Tally_b(t) \,\&\, \forall t'(Tally_b(t') \,\&\, t' \subseteq_p w \rightarrow t' \subseteq_p t).$$

Let us say, further, when a b-tally t is <u>longer than any b-tally in</u> x:

$$Max^+T_b(t,x) \equiv MaxT_b(t,x) \,\&\, \neg t \subseteq_p x.$$

We then define when a string u is a <u>preframe indexed by</u> t:

$$Pref(u,t) \equiv \exists y \subseteq_p u\, (aya=u \,\&\, Max^+T_b(t,u));$$

when $t_1ut_2$ is (the) <u>first frame</u> in the string x, $Firstf(x,t_1,u,t_2)$:

$$Pref(u,t_1) \,\&\, Tally_b(t_2) \,\&\, ((t_1=t_2 \,\&\, t_1ut_2=x) \vee (t_1<t_2 \,\&\, (t_1ut_2a)Bx));$$

when $t_1ut_2$ is (the) <u>last frame</u> in x, $Lastf(x,t_1,u,t_2)$:

$$Pref(u,t_1) \,\&\, t_1=t_2 \,\&\, (t_1ut_2=x \vee \exists w\, (wat_1ut_2=x \,\&\, Max^+T_b(t_1,w)));$$

and when $t_1ut_2$ is an <u>intermediate frame</u> in x immediately following an initial segment w of x, $Intf(x,w,t_1,u,t_2)$:

$$Pref(u,t_1) \,\&\, Tally_b(t_2) \,\&\, t_1<t_2 \,\&\, \exists w_1(wat_1ut_2aw_1=x) \,\&\, Max^+T_b(t_1,w).$$

Then we define when a string u is <u>$t_1,t_2$-framed</u> in x:

$Fr(x,t_1,u,t_2) \equiv Firstf(x,t_1,u,t_2)\ \ v\ \ \exists w\ Intf(x,w,t_1,u,t_2)\ \ v\ \ Lastf(x,t_1,u,t_2)$,

We say that $t_1$ is the <u>initial</u>, and $t_2$ <u>terminal tally marker</u> in the frame.

Next we define "t <u>envelops</u> x", Env(t,x), to be the conjunction of the following five conditions:

(a) $MaxT_b(t,x)$ — "t is a longest b-tally in x",

(b) $\exists u \subseteq_p x\ \exists t_1,t_2\ Firstf(x,t_1,u,t_2)$ — "x has a first frame",

(c) $\exists u \subseteq_p x\ Lastf(x,t,u,t)$ — "x has a last frame with t as its initial and terminal marker"

(d) $\forall u \subseteq_p x\ \forall t_1,t_2,t_3,t_4\ (Fr(x,t_1,u,t_2)\ \&\ Fr(x,t_3,u,t_4)\ \rightarrow\ t_1=t_3)$

"different initial tally markers frame distinct strings",

(e) $\forall u_1,u_2 \subseteq_p x\ \forall t',t_1,t_2\ (Fr(x,t',u_1,t_1)\ \&\ Fr(x,t',u_2,t_2)\ \rightarrow\ u_1=u_2)$

"distinct strings are framed by different initial tally markers"

Now we say x <u>is a set code</u> if x is aa or else x is enveloped by some b-tally:

$$Set(x)\ \equiv\ x=aa\ \ v\ \ \exists t \subseteq_p x\ Env(t,x).$$

Finally, we say that a string y <u>is a member of the set coded by</u> string x if x is enveloped by some b-tally t and the juxtaposition of the string y with single tokens of digit a is framed in x:

$$y\ \varepsilon\ x\ \ \equiv\ \ \exists t \subseteq_p x\ (Env(t,x)\ \&\ \exists u \subseteq_p x\ \exists t_1,t_2(Fr(x,t_1,u,t_2)\ \&\ u=aya)).$$

We can then establish:

SINGLETON LEMMA. For any string form $I \subseteq I_0$ there is a string form $J \subseteq I$ such that

$$QT^+ \vdash \forall x \in J\ \forall u,t_1,t_2\ (Set(x)\ \&\ Firstf(x,t_1,aua,t_2)\ \&\ x=t_1auat_2 \rightarrow$$

$$\rightarrow \forall w\ (w\ \varepsilon\ x \leftrightarrow w=u)).$$

(See [2], (5.21).)

APPENDING LEMMA. For any string form $I \subseteq I_0$ there is a string form $J \subseteq I$ such that

$$QT^+ \vdash \forall x,y \in J\ \forall t,t_2,t_3(Env(t_2,x)\ \&\ Env(t,y)\ \&\ (t_3a)By\ \&\ Tally_b(t_3)\ \&\ t_2<t_3\ \&$$

$$\&\ \neg\exists u(u\ \varepsilon\ x\ \&\ u\ \varepsilon\ y) \rightarrow \exists z \in J\ (Env(t,z)\ \&\ \forall u(u\ \varepsilon\ z \leftrightarrow u\ \varepsilon\ x \vee u\ \varepsilon\ y)).$$

(See [2], (5.46).)

We then derive:

DOUBLETON LEMMA. For any string form $I \subseteq I_0$ there is a string form $J \subseteq I$ such that

$$QT^+ \vdash \forall x \in J\ \forall t_1,t_2,t_3,u,v(Pref(aua,t_1)\ \&\ Pref(ava,t_2)\ \&\ t_1<t_2\ \&\ t_2=t_3\ \&\ u \neq v\ \&$$

$$\&\ x=t_1auat_2avat_3 \rightarrow Set(x)\ \&\ \forall w(w\ \varepsilon\ x \leftrightarrow (w=u \vee w=v)).$$

(See [2], (5.58).)

Let $\quad \text{MinMax}^{+}T_b(t,u) \equiv \text{Max}^{+}T_b(t,u) \ \&\ \forall t'(\text{Max}^{+}T_b(t',u) \rightarrow t \leq t')$.

In that case we say that t is a shortest non-occurrent b-tally in string u.

We then have:

SHORTEST NON-OCCURRENT TALLY LEMMA. For any string form $I \subseteq I_0$ there is a string form $J \subseteq I$ such that

$$QT^{+} \vdash \forall x \in J\ \exists! t \in J\ \text{MinMax}^{+}T_b(t,x).$$

STRING RECURSION THEOREM. Let $F_1(y,z,u)$ and $F_2(y,z,u)$ be formulae, and let $I \subseteq I^{\Diamond}$ closed under * and downward closed under $\subseteq_p$. Suppose that

$$QT^{+} \vdash I(p) \ \&\ I(q),$$

$$QT^{+} \vdash \forall y,z \in I\ \exists! u \in I\ F_1(y,z,u),$$

and $\quad QT^{+} \vdash \forall y,z \in I\ \exists! u \in I\ F_2(y,z,u)$.

Then there is a formula $H(y,z)$ and a string form $J \subseteq I$ such that

(i) $QT^+ \vdash \forall y \in J\ \exists! z \in I\ H(y,z)$,

(iia) $QT^+ \vdash \forall y \in I\ (H(a,y) \leftrightarrow y=p)$,

(iib) $QT^+ \vdash \forall y \in I\ (H(b,y) \leftrightarrow y=q)$,

(iiia) $QT^+ \vdash \forall y \in J\ \forall u,z \in I\ (H(y,u) \to (H(y{*}a,z) \leftrightarrow F_1(y,u,z)))$,

and (iiib) $QT^+ \vdash \forall y \in J\ \forall u,z \in I\ (H(y,u) \to (H(y{*}b,z) \leftrightarrow F_2(y,u,z)))$.

<u>Proof</u>: Let Comp(u,m) abbreviate

$$\begin{aligned}
&Set(u)\ \&\ (a \le m \to \exists v \subseteq_p u\ (Pair[a,p,v]\ \&\ v\ \varepsilon\ u))\ \& \\
&\quad \&\ (b \le m \to \exists v \subseteq_p u\ (Pair[b,q,v]\ \&\ v\ \varepsilon\ u))\ \& \\
&\ \&\ \forall z<m\ \forall u_1,u_2,v_1\ (Pair[z,u_1,v_1]\ \&\ v_1\ \varepsilon\ u\ \&\ F_1(z,u_1,u_2) \to \\
&\qquad\qquad \to \exists v_2 \subseteq_p u\ (Pair[z{*}a,u_2,v_2]\ \&\ v_2\ \varepsilon\ u))\ \& \\
&\ \&\ \forall z<m\ \forall u_1,u_2,v_1\ (Pair[z,u_1,v_1]\ \&\ v_1\ \varepsilon\ u\ \&\ F_2(z,u_1,u_2) \to \\
&\qquad\qquad \to \exists v_2 \subseteq_p u\ (Pair[z{*}b,u_2,v_2]\ \&\ v_2\ \varepsilon\ u))\ \& \\
&\&\ \forall z,u_1,u_2,v_1,v_2\ (Pair[z,u_1,v_1]\ \&\ Pair[z,u_2,v_2]\ \&\ v_1\ \varepsilon\ u\ \&\ v_2\ \varepsilon\ u \to \\
&\qquad\qquad \to u_1=u_2\ \&\ v_1=v_2).
\end{aligned}$$

Let (C1)-(C6) be the successive conjuncts that make up Comp(u,m,x). Then (C4) and (C5) express the usual conditions that a sequence code u should satisfy to represent the course of a recursion. The last clause, (C6), is a uniqueness condition. Then Comp(u,m) means, roughly, that u is a set code for a computation determined by p, q, $F_1$,$F_2$, in at least m steps where the length indices m are strings ordered by the tree-like ordering ≤.

Let MinComp(u,m) abbreviate

$$\mathrm{Comp}(u,m)\ \&\ \forall u'\ (\mathrm{Comp}(u',m) \rightarrow \forall y\ (y\ \varepsilon\ u \rightarrow y\ \varepsilon\ u'))\ \&$$
$$\&\ \forall z,v,w\ (\mathrm{Pair}[z,v,w]\ \&\ w\ \varepsilon\ u \rightarrow (m=a\ \&\ z=a)\ v\ (m=b\ \&\ z=b)\ v$$
$$v\ \exists n<m\ (z\leq na\ v\ z\leq nb)).$$

Let J(m) abbreviate

$$I(m)\ \&\ \exists!y\in I\ \exists u\in I\ \exists w\subseteq_p u\ (\mathrm{MinComp}(u,m)\ \&\ \mathrm{Pair}[m,y,w]\ \&\ w\ \varepsilon\ u).$$

Finally, let H(m,y) abbreviate

$$\exists u,w\ (\mathrm{MinComp}(u,m)\ \&\ \mathrm{Pair}[m,y,w]\ \&\ w\ \varepsilon\ u).$$

Let (C1)-(C6) be the successive conjuncts that make up Comp(u,m). Then (C4) and (C5) express the usual conditions that a sequence code u should satisfy to represent the course of a recursion. The last clause, (C6), is a uniqueness condition.

The proof consists of ten claims.

Claim 1: $QT^+ \vdash J(a)$.

By the principal hypothesis, $QT^+ \vdash I(p)$.

By the Pairing Lemma, $QT^+ \vdash \exists! w \in I\ Pair[a,p,w]$.

By the Shortest Non-Occurrent Tally Lemma, $M \models \exists! t \in I\ MinMax^+T_b(t,awa)$.

$\Rightarrow M \models Max^+T_b(t,awa)$.

Let $u=tawat$.

Then $M \models I(u)$.

$\Rightarrow M \models Firstf(u,t,awa,t)\ \&\ Lastf(u,t,awa,t)$,

$\Rightarrow$ by the Singleton Lemma, $M \models Set(u)\ \&\ \forall z\ (z\ \varepsilon\ u \leftrightarrow z=w)$,

$\Rightarrow M \models Set(u)\ \&\ Pair[a,p,w]\ \&\ w\ \varepsilon\ u$,

which suffices to establish parts (C1) and (C2) of $M \models Comp(u,a)$.

Since $QT^+ \vdash \neg(b \leq a)$ and $QT^+ \vdash \forall z\ \neg(z<a)$, parts (C3)-(C5) hold trivially.

For (C6), assume that

$M \models Pair[z,u_1,v_1]\ \&\ Pair[z,u_2,v_2]\ \&\ v_1\ \varepsilon\ u\ \&\ v_2\ \varepsilon\ u$.

$\Rightarrow$ by choice of u, $M \models v_1=v_2$,

$\Rightarrow$ since $M \models v_1 \subseteq_p u\ \&\ I(u)$, $M \models I(v_1)$,

$\Rightarrow M \models Pair[z,u_1,v_1]\ \&\ Pair[z,u_2,v_1]$,

$\Rightarrow$ by the Pairing Lemma, $M \models u_1=u_2$,

$\Rightarrow M \models u_1=u_2\ \&\ v_1=v_2$, as required.

This completes the argument that $M \models Comp(u,a)$. We now move on to show that $M \models MinComp(u,a)$.

Assume now that $M \models Comp(v,a)$.

Then $M \models \exists w_1 \subseteq_p v\ (Pair[a,p,w_1]\ \&\ w_1\ \varepsilon\ v)$ .

$\Rightarrow M \models Pair[a,p,w]\ \&\ Pair[a,p,w_1]$,

$\Rightarrow$ by the Pairing Lemma, $M \models w=w_1$,

$\Rightarrow M \models w\ \varepsilon\ v$.

Assume now that $M \models y\ \varepsilon\ u$.

$\Rightarrow$ from $M \models \forall z\ (z\ \varepsilon\ u \leftrightarrow z=w)$, $M \models y=w$,

$\Rightarrow M \models y\ \varepsilon\ v$.

Thus we proved that $M \models Comp(v,a) \rightarrow \forall y(y\ \varepsilon\ u \rightarrow y\ \varepsilon\ v)$.

To complete the argument that $M \models MinComp(u,a)$, assume that

$$M \models Pair[z_1,v_1,w_1]\ \&\ w_1\ \varepsilon\ u.$$

$\Rightarrow$ from $M \models \forall z\ (z\ \varepsilon\ u \leftrightarrow z=w)$, $M \models w_1=w$,

$\Rightarrow M \models Pair[a,p,w]\ \&\ Pair[z_1,v_1,w]$,

$\Rightarrow$ by the Pairing Lemma, $M \models z_1=a\ \&\ v_1=p$.

Therefore we also have

$$M \models \forall z_1,v_1,w_1\ (Pair[z_1,v_1,w_1]\ \&\ w_1\ \varepsilon\ u\ \rightarrow (a=a\ \&\ z_1=a)\ v\ (a=b\ \&\ z_1=b)\ v$$
$$v\ \exists n<a\ ((z_1 \le na\ v\ z_1 \le nb))).$$

So we finally have that $M \models MinComp(u,a)$.

In fact, we obtained

$M \models \exists! y \in I\ \exists u \in I\ \exists w \subseteq_p u\ (MinComp(u,a)\ \&\ Pair[a,y,w]\ \&\ w\ \varepsilon\ u)$.

So $M \vDash J(a)$.

Claim 2: $QT^+ \vdash J(b)$.

From the proof of $QT^+ \vdash J(a)$ we have that

$$M \vDash \exists! w_1 \in I\ (Pair[a,p,w_1]\ \&\ \exists! t_1 \in I\ MinMax^+T_b(t_1,aw_1a)).$$

Arguing exactly analogously, we obtain that

$$M \vDash \exists! w_2 \in I\ (Pair[b,q,w_2]\ \&\ \exists! t_2 \in I\ MinMax^+T_b(t_2,aw_2a)).$$

$\Rightarrow$ since $QT^+ \vdash a \neq b$, from $M \vDash Pair[a,p,w_1]\ \&\ Pair[b,q,w_2]$, by the Pairing Lemma, $M \vDash w_1 \neq w_2$.

Let $u' = t_1aw_1a(t_1t_2)aw_2a(t_1t_2)$. Then $M \vDash I(u')$.

$\Rightarrow$ by the proof of Doubleton Lemma,

$$M \vDash Env(t_1t_2,u')\ \&\ \forall w\ (w\ \varepsilon\ u' \leftrightarrow w=w_1 \vee w=w_2).$$

On the other hand, by the principal hypothesis,

$$M \vDash \exists! u_3 \in I\ F_1(a,u_1,u_3) \text{ and } M \vDash \exists! u_4 \in I\ F_2(a,u_1,u_4).$$

Just as above, we then obtain

$$M \vDash \exists! w_3 \in I\ (Pair[aa,u_3,w_3]\ \&\ \exists! t_3 \in I\ MinMax^+T_b(t_3,aw_3a))$$

and $$M \vDash \exists! w_4 \in I\ (Pair[ab,u_4,w_4]\ \&\ \exists! t_4 \in I\ MinMax^+T_b(t_4,aw_4a)).$$

Then, just as above, we again have that $M \vDash w_3 \neq w_4$, and, further, that

$$M \vDash w_1 \neq w_3\ \&\ w_1 \neq w_4\ \&\ w_2 \neq w_3\ \&\ w_2 \neq w_4.$$

Letting $u'' = (t_1t_2t_3)aw_3a(t_1t_2t_3t_4)aw_4a(t_1t_2t_3t_4)$, we likewise have that $M \vDash I(u'')$ and $M \vDash Env(t_1t_2t_3t_4,u'')\ \&\ \forall w\ (w\ \varepsilon\ u'' \leftrightarrow w=w_3 \vee w=w_4)$.

Since

$M \vDash Env(t_1t_2,u') \;\&\; Env(t_1t_2t_3t_4,u'') \;\&\; (t_1t_2t_3a)Bu'' \;\&$

$$\& \; Tally_b(t_1t_2t_3) \;\&\; \neg\exists w(w \,\varepsilon\, u' \;\&\; w \,\varepsilon\, u''),$$

it follows by the proof of Appending Lemma, that for

$$u=t_1aw_1a(t_1t_2)aw_2a(t_1t_2t_3)aw_3a(t_1t_2t_3t_4)aw_4a(t_1t_2t_3t_4),$$

we have $M \vDash Env(t_1t_2t_3t_4,u) \;\&\; \forall w\,(w \,\varepsilon\, u \leftrightarrow w \,\varepsilon\, u' \vee w \,\varepsilon\, u'')$.

Hence

$$M \vDash Set(u) \;\&\; \forall w\,(w \,\varepsilon\, u \leftrightarrow w=w_1 \vee w=w_2 \vee w=w_3 \vee w=w_4).$$

So (C1) holds.

Now, we have that

$$M \vDash \exists w_1 \subseteq_p u\,(Pair[a,p,w_1] \;\&\; w_1 \,\varepsilon\, u)$$

and $M \vDash \exists w_2 \subseteq_p u\,(Pair[b,q,w_2] \;\&\; w_2 \,\varepsilon\, u)$.

Since $QT^+ \vdash a \leq b$, this suffices to establish (C2) and (C3) of $M \vDash Comp(u,b)$.

Since $QT^+ \vdash \forall z(z<b \rightarrow z=a)$, and $M \vDash w_3 \,\varepsilon\, u \;\&\; w_4 \,\varepsilon\, u$, we have from the choices of $w_3$ and $w_4$, that (C4) and (C5) of $M \vDash Comp(u,b)$ also hold.

For (C6), assume that

$$M \vDash Pair[z,s_1,v_1] \;\&\; Pair[z,s_2,v_2] \;\&\; v_1 \,\varepsilon\, u \;\&\; v_2 \,\varepsilon\, u.$$

$\Rightarrow$ $M \vDash (v_1=w_1 \vee v_1=w_2 \vee v_1=w_3 \vee v_1=w_4) \;\&$

$$\& \; (v_2=w_1 \vee v_2=w_2 \vee v_2=w_3 \vee v_2=w_4).$$

Suppose, for a reductio that $M \vDash v_1 \neq v_2$, say $M \vDash v_1=w_1 \;\&\; v_2=w_2$.

$\Rightarrow$ $M \vDash Pair[z,s_1,w_1] \;\&\; Pair[z,s_2,w_2]$.

But we have that $M \vDash Pair[a,p,w_1] \;\&\; Pair[b,q,w_2]$.

⇒ by the Pairing Lemma, $M \vDash z=a \;\&\; z=b$, a contradiction.

Similarly for the other choice of for $v_1$ and $v_2$ from u.

Therefore $M \vDash v_1=v_2$, and (C6) of $M \vDash \text{Comp}(u,b)$ also holds.

Assume now that $M \vDash \text{Comp}(v,b)$.

Then

$$M \vDash \exists p_1 \subseteq_p v\, (\text{Pair}[a,p,p_1] \;\&\; p_1 \,\varepsilon\, v), \quad \text{and}$$

$$M \vDash \exists p_2 \subseteq_p v\, (\text{Pair}[b,q,p_2] \;\&\; p_2 \,\varepsilon\, v), \quad \text{and}$$

$$M \vDash \exists v_3,p_3 \subseteq_p v\, (\, F_1(a,v_1,v_3) \;\&\; I(v_3) \;\&\; \text{Pair}[aa,v_3,p_3] \;\&\; p_3 \,\varepsilon\, v), \quad \text{and}$$

$$M \vDash \exists v_4,p_4 \subseteq_p v\, (\, F_2(a,v_1,v_4) \;\&\; I(v_4) \;\&\; \text{Pair}[ab,v_4,p_4] \;\&\; p_4 \,\varepsilon\, v).$$

From the principal hypothesis, it follows that

$$M \vDash v_3=u_3 \;\&\; v_4=u_4.$$

⇒ $M \vDash \text{Pair}[a,p,w_1] \;\&\; \text{Pair}[a,p,p_1]$ and $M \vDash \text{Pair}[b,q,w_2] \;\&\; \text{Pair}[b,q,p_2]$

and $M \vDash \text{Pair}[aa,u_3,w_3] \;\&\; \text{Pair}[aa,u_3,p_3]$ and $M \vDash \text{Pair}[ab,u_4,w_4] \;\&\; \text{Pair}[ab,u_4,p_4]$,

⇒ by the Pairing Lemma, $M \vDash w_1=p_1 \;\&\; w_2=p_2 \;\&\; w_3=p_3 \;\&\; w_4=p_4$,

⇒ from $M \vDash w_1 \,\varepsilon\, u \;\&\; w_2 \,\varepsilon\, u \;\&\; w_3 \,\varepsilon\, u \;\&\; w_4 \,\varepsilon\, u$,

$$M \vDash w_1 \,\varepsilon\, v \;\&\; w_2 \,\varepsilon\, v \;\&\; w_3 \,\varepsilon\, v \;\&\; w_4 \,\varepsilon\, v.$$

Assume now that $M \vDash y \,\varepsilon\, u$.

⇒ $M \vDash y=w_1 \vee y=w_2 \vee y=w_3 \vee y=w_4$,

⇒ $M \vDash y \,\varepsilon\, v$.

Thus $M \vDash \text{Comp}(v,b) \rightarrow \forall y(y \,\varepsilon\, u \rightarrow y \,\varepsilon\, v)$.

Finally, assume that

$$M \vDash \text{Pair}[z,s,w] \;\&\; w \,\varepsilon\, u.$$

⇒ M ⊨ $w=w_1$ v $w=w_2$ v $w=w_3$ v $w=w_4$,

⇒ M ⊨ Pair[z,s,$w_1$] v Pair[z,s,$w_2$] v Pair[z,s,$w_3$] v Pair[z,s,$w_4$].

But we have

M ⊨ Pair[a,p,$w_1$] & Pair[b,q,$w_2$] & Pair[aa,$u_3$,$w_3$] & Pair[ab,$u_4$,$w_4$].

⇒ by the Pairing Lemma, M ⊨ z=a v z=b v z=aa v z=ab.

We have that M ⊨ a<b.

Hence, from M ⊨ m=b, since M ⊨ a≤ab & aa≤aa & ab≤ab, it follows that

M ⊨ (m=b & z=b) v ∃n<b (z≤na v z≤nb)

whence M ⊨ (m=a & z=a) & (m=b & z=b) v ∃n<b (z≤na v z≤nb),

as required.

This completes the argument that M ⊨ MinComp(u,b).

<u>Claim 3</u>: $QT^+$ ⊢ ∀x (J(x) → J(Sx)).

Assume that M ⊨ J(m).

⇒ M ⊨ I(m),

⇒ since I is a string form, M ⊨ I(Sm).

We need to show that

M ⊨ ∃!y∈I ∃u∈I ∃w$\subseteq_p$u (MinComp(u,Sm) & Pair[Sm,y,w] & w ε u).

If M ⊨ Sm=b, what we need was proved in Claim 2.

So we may assume M ⊨ ¬Sm=b. Then M ⊨ ¬m=a.

From the hypothesis M ⊨ J(m) we have that

M ⊨ ∃!y∈I ∃u∈I ∃w$\subseteq_p$u (MinComp(u,m) & Pair[m,y,w] & w ε u).

Let $u_0$ be a u in M such that $M \vDash I(u) \ \&\ Set(u) \ \&\ MinComp(u,m,x)$.

Let $y_0$ be the unique y in M such that

$$M \vDash \exists w{\subseteq_p}u_0\ (Pair[m,y,w] \ \&\ w \ \varepsilon\ u_0).$$

$\Rightarrow$ since $M \vDash I(m) \ \&\ I(y_0)$, by the Pairing Lemma, $M \vDash \exists! w{\subseteq_p}u_0\ Pair[m,y,w]$.

Let $w_0$ be the unique such w in M.

From $M \vDash MinComp(u_0,m)$, $M \vDash Comp(u_0,m)$.

$\Rightarrow$ $M \vDash Set(u_0)$,

$\Rightarrow$ since $M \vDash w_0 \ \varepsilon\ u_0$, $M \vDash \exists t{\subseteq_p}u_0\ Env(t,u_0)$.

Here t uniquely depends on $u_0$. (See [2], $(4.24^b)$.) Since the string form I is downward closed w.r.t. to $\subseteq_p$, from $M \vDash I(u_0)$ we have that $M \vDash I(t)$.

From the principal hypothesis of the Theorem we have that

(†) $M \vDash \exists! v_1{\in}I\ F_1(m,y_0,v_1)$ and $M \vDash \exists! v_2{\in}I\ F_2(m,y_0,v_2)$.

$\Rightarrow$ by the Pairing Lemma,

$$M \vDash \exists! w_1{\in}I\ (Pair[ma,v_1,w_1] \ \&\ \exists! t_1{\in}I\ MinMax^+T_b(t_1,aw_1a)) \text{ and}$$

$$M \vDash \exists! w_2{\in}I\ (Pair[mb,v_2,w_2] \ \&\ \exists! t_2{\in}I\ MinMax^+T_b(t_2,aw_2a)).$$

Then, analogously to the proof of $QT^+ \vdash J(b)$ above, we obtain, for

$$u' = tt_1aw_1a(tt_1t_2)aw_2a(tt_1t_2),$$

that $M \vDash I(u') \ \&\ Env(tt_1t_2,u') \ \&\ \forall w(w \ \varepsilon\ u' \leftrightarrow w{=}w_1 \vee w{=}w_2)$.

From $M \vDash MinComp(u_0,m)$, we readily verify that

$$M \vDash \neg(w_1 \ \varepsilon\ u_0) \ \&\ \neg(w_2 \ \varepsilon\ u_0).$$

Then, since

$M \vDash Env(t,u_0) \ \&\ Env(tt_1t_2,u') \ \&\ (tt_1t_2a)Bu' \ \&\ Tally_b(tt_1t_2) \ \&$

$\& \neg\exists w(w \,\varepsilon\, u_0 \,\&\, w \,\varepsilon\, u'),$

by the proof of Appending Lemma, for

$$u=u_0t_1aw_1a(tt_1t_2)aw_2a(tt_1t_2t_3),$$

we have $M \models Env(tt_1t_2) \,\&\, \forall w\, (w \,\varepsilon\, u \leftrightarrow w \,\varepsilon\, u_0 \vee w \,\varepsilon\, u')$.

Hence

$$M \models Set(u) \,\&\, \forall w\, (w \,\varepsilon\, u \leftrightarrow (w \,\varepsilon\, u_0 \vee w=w_1 \vee w=w_2)).$$

So (C1) holds.

Note that, since the string form I is closed under *, from

$$M \models I(u_0) \,\&\, I(t_1) \,\&\, I(w_1) \,\&\, I(t) \,\&\, I(t_2) \,\&\, I(w_2)$$

we have $M \models I(u)$.

We now proceed to argue that $M \models Comp(u,Sm)$.

It is straightforward to verify from $M \models Comp(u_0,m)$ and the choice of u that

$M \models \exists q_1 \subseteq_p u\, (Pair[a,p,q_1] \,\&\, q_1 \,\varepsilon\, u)$ and $M \models \exists q_2 \subseteq_p u\, (Pair[b,q,q_2] \,\&\, q_2 \,\varepsilon\, u)$

so that (C2) and (C3) of $M \models Comp(u,Sm)$ both hold.

For (C4), let $M \models z<Sm \,\&\, Pair[z,u_1,v_3] \,\&\, v_3 \,\varepsilon\, u \,\&\, F_1(z,u_1,u_2)$ where

$M \models u_1,u_2,v_3 \subseteq_p u$.

We need to show that $M \models \exists v \subseteq_p u\, (Pair[z^*a,u_2,v] \,\&\, v \,\varepsilon\, u)$.

From $M \models z<Sm$, $M \models z<m \vee z=m$.

Suppose $M \models z<m$.

We have that $M \models v_3 \,\varepsilon\, u$.

Using the Pairing Lemma and the definition of $<$, we verify that

$$M \models v_3 \neq w_1 \text{ and } M \models v_3 \neq w_2.$$

⇒ from $M \vDash v_3 \varepsilon u$, $M \vDash v_3 \varepsilon u_0$,

⇒ $M \vDash u_1 \subseteq_p u_0$.

Then, from $M \vDash Pair[z,u_1,v_3] \& v_3 \varepsilon u_0 \& F_1(z,u_1,u_2)$ and $M \vDash Comp(u_0,m)$,

we have that $M \vDash \exists v \subseteq_p u_0 (Pair[z^*a,u_2,v] \& v \varepsilon u_0)$,

whence $M \vDash \exists v \subseteq_p u (Pair[z^*a,u_2,v] \& v \varepsilon u)$, as required.

Suppose $M \vDash z=m$.

Again, we are assuming that $M \vDash Pair[z,u_1,v_3] \& v_3 \varepsilon u \& F_1(z,u_1,u_2)$ where $M \vDash u_1,u_2,v_3 \subseteq_p u$. Hence $M \vDash I(u_2)$.

Just as above, $M \vDash v_3 \varepsilon u_0$.

On the other hand, we also have that $M \vDash \exists w \subseteq_p u (Pair[m,y_0,w] \& w \varepsilon u_0)$,

⇒ from $M \vDash Pair[m,u_1,v_3] \& v_3 \varepsilon u_0$ and clause (C6) of $M \vDash Comp(u_0,m)$,

$$M \vDash y_0=u_1 \& v_3=w,$$

⇒ from $M \vDash F_1(m,u_1,u_2) \& I(u_2)$ and (†), $M \vDash v_1=u_2$.

But then, from $M \vDash Pair[ma,v_1,w_1]$, we have

$$M \vDash Pair[ma,u_2,w_1] \& w_1 \varepsilon u,$$

where $M \vDash w_1 \subseteq_p u$, as required.

Hence (C4) of $M \vDash Comp(u,Sm)$ also holds.

For (C5), we argue in exactly the same way, except that references to $F_1$, $z^*a$ and $w_1$ are replaced by $F_2$, $z^*b$ and $w_2$.

Condition (C6) is verified using the corresponding condition from $M \vDash Comp(u_0,m)$ and the Pairing Lemma.

We now proceed to show that in fact $M \vDash MinComp(u,Sm)$.

Suppose that $M \models Comp(v',Sm)$.

First, we want to show that $M \models \forall y(y \varepsilon v' \rightarrow y \varepsilon u)$.

From $M \models Comp(v',Sm)$ we have that $M \models Comp(v',m)$.

From the hypothesis $M \models J(m)$ we have that $M \models MinComp(u_0,m)$.

$\Rightarrow M \models \forall y(y \varepsilon u_0 \rightarrow y \varepsilon v')$.

From $M \models J(m)$ we also have that $M \models Pair[m,y_0,w_0] \ \& \ w_0 \varepsilon u_0$.

$\Rightarrow M \models w_0 \varepsilon v'$.

But then, since, by (†), $M \models F_1(m,y_0,v_1) \ \& \ F_2(m,y_0,v_2)$,

we have, from (C4) and (C5) of $M \models Comp(v',Sm)$ that

$$M \models w_1 \varepsilon v' \ \& \ w_2 \varepsilon v'$$

where $M \models Pair[m*a,v_1,w_1] \ \& \ Pair[m*b,v_2,w_2]$.

So we have that $M \models \forall y(y \varepsilon u_0 \rightarrow y \varepsilon v') \ \& \ w_1 \varepsilon v' \ \& \ w_2 \varepsilon v'$.

But then from, the choice of u, it follows that

$$M \models \forall y(y \varepsilon u \rightarrow y \varepsilon v'),$$

as required.

Suppose now that $M \models Pair[z,v,w] \ \& \ w \varepsilon u$.

$\Rightarrow M \models w \varepsilon u_0 \ v \ w{=}w_1 \ v \ w{=}w_2$.

If $M \models w \varepsilon u_0$, then from $M \models MinComp(u_0,m)$, we have that

$$M \models (m{=}a \ \& \ z{=}a) \ \& \ (m{=}b \ \& \ z{=}b) \ v \ \exists n{<}m \ (z{\leq}na \ v \ z{\leq}nb).$$

But then $M \models (m{=}a \ \& \ z{=}a) \ \& \ (m{=}b \ \& \ z{=}b) \ v \ \exists n{<}Sm \ (z{\leq}na \ v \ z{\leq}nb)$.

If $M \models w{=}w_1$, we have that $M \models z{=}ma$, whence

$$M \models \exists n{<}Sm \ z{\leq}na.$$

Hence $M \vDash (m=a \,\&\, z=a) \,\&\, (m=b \,\&\, z=b) \vee \exists n<Sm\, (z \leq na \vee z \leq nb)$,

as required.

An analogous argument applies if $M \vDash w=w_2$.

This suffices to establish $M \vDash MinComp(u,Sm)$.

Now, we have that

$M \vDash \exists! w_2 \in I\, (MinComp(u,Sm) \,\&\, Pair[Sm,v_2,w_2] \,\&\, w_2\, \varepsilon\, u)$.

Suppose that $M \vDash MinComp(u,Sm) \,\&\, Pair[Sm,y,w_2] \,\&\, w_2\, \varepsilon\, u$ where $M \vDash I(y)$.

$\Rightarrow$ from $M \vDash Pair[Sm,v_2,w_2] \,\&\, I(w_2)$, we have, by the Pairing Lemma, that

$$M \vDash v_2=y.$$

So we have actually established that

$M \vDash \exists! y \in I\, \exists u \in I\, \exists w \subseteq_p u\, (MinComp(u,Sm) \,\&\, Pair[Sm,y,w] \,\&\, w\, \varepsilon\, u)$,

and hence that $M \vDash J(Sm)$.

This completes the proof of Claim 3.

<u>Claim 4</u>: $QT^+ \vdash \forall x\, (J(x) \rightarrow J(x*a))$.

Exactly analogous to the proof of Claim 3.

Claims 1-4 establish that J is a string form.

<u>Claim 5</u>: $QT^+ \vdash \forall y \in I\, (H(a,y) \leftrightarrow y=p)$.

Let $M \vDash I(y)$.

Assume $M \vDash y=p$.

As shown in the proof of Claim 1, in M there is a u, namely, u=tawat, such that

$M \models Pair[a,p,w] \ \& \ w \ \varepsilon \ u.$

Then, again as shown in the proof of Claim 1, we have that

$M \models MinComp(u,a) \ \& \ Pair[a,p,w] \ \& \ w \ \varepsilon \ u),$

whence $M \models H(a,p)$, so $M \models H(a,y)$.

Thus, $M \models y=p \rightarrow H(x,a,y)$.

Conversely, let $M \models H(a,y)$.

$\Rightarrow$ by definition of H, $M \models \exists u,w \ (MinComp(u,a) \ \& \ Pair[a,y,w] \ \& \ w \ \varepsilon \ u),$

$\Rightarrow M \models Comp(u,a),$

$\Rightarrow$ from (C2), $M \models \exists v \ (Pair[a,p,v] \ \& \ v \ \varepsilon \ u),$

$\Rightarrow$ from $M \models Pair[a,y,w] \ \& \ Pair[a,p,v] \ \& \ w \ \varepsilon \ u \ \& \ v \ \varepsilon \ u$, and clause (C6) of

$M \models Comp(u,a)$, $M \models y=p$.

Hence also $M \models H(a,y) \rightarrow y=p)$.

This completes the proof of Claim 5.

<u>Claim 6</u>: $QT^+ \vdash \forall y \in I \ (H(b,y) \leftrightarrow y=q)$.

Let $M \models I(y)$.

Assume $M \models y=q$.

We follow the proof of Claim 2 to obtain a u in M such that

$M \models Pair[b,u_2,w] \ \& \ w \ \varepsilon \ u,$

where $M \models MinComp(u,b) \ \& \ Pair[b,u_2,w] \ \& \ w \ \varepsilon \ u).$

Then $M \models H(b,y)$.

This shows that $M \models y=q \rightarrow H(b,y)$.

To establish the converse, that $M \models H(b,y) \rightarrow y=q$, we argue analogously to the proof in Claim 5.

Claim 7: $QT^+ \vdash \forall y \in J\ \forall v,z \in I\ (H(y,v) \rightarrow (F_1(y,v,z) \rightarrow H(y*a,z)))$.

Let $M \models J(y)$ and $M \models I(v)\ \&\ I(z)$.

Suppose $M \models H(y,v)\ \&\ F_1(y,v,z)$.

$\Rightarrow$ from $M \models J(y)$,

$$M \models \exists u_0 \in I\ \exists w \subseteq_p u\ (\mathrm{MinComp}(u_0,y)\ \&\ \mathrm{Pair}[y,v,w]\ \&\ w\ \varepsilon\ u_0).$$

We then obtain, exactly analogously to the proof of Claim 3, a u in M such that

$$M \models \exists w_1\ (\mathrm{MinComp}(u,y*a)\ \&\ \mathrm{Pair}[y*a,z,w_1]\ \&\ w_1\ \varepsilon\ u),$$

whence $M \models H(y*a,z)$ follows.

This completes the argument for Claim 7.

Claim 8: $QT^+ \vdash \forall y \in J\ \forall v,z \in I\ (H(y,v)\ \&\ H(y*a,z) \rightarrow F_1(y,v,z))$.

Assume that $M \models H(y,v)\ \&\ H(y*a,z)$ where $M \models J(y)$ and $M \models I(v)\ \&\ I(z)$.

From the hypothesis $M \models H(y,v)$ we have that

$$M \models \exists u_0,w_0\ (\mathrm{MinComp}(u_0,y)\ \&\ \mathrm{Pair}[y,v,w_0]\ \&\ w_0\ \varepsilon\ u_0).$$

From the principal hypothesis of the Theorem we have

$QT^+ \vdash \exists! z' \in I\ F_1(y,v,z')$.

We then obtain, exactly analogously to the proof of Claim 3, a u in M such that

$M \vDash \exists w_1\ (\mathrm{MinComp}(u,y*a)\ \&\ \mathrm{Pair}[y*a,z',w_1]\ \&\ w_1\ \varepsilon\ u)$.

On the other hand, from the hypothesis $M \vDash H(y*a,z)$, we have that

$M \vDash \exists u',w'\ (\mathrm{MinComp}(u',y*a)\ \&\ \mathrm{Pair}[y*a,z',w']\ \&\ w'\ \varepsilon\ u')$.

Now, we have that, in general

$QT^+ \vdash \mathrm{MinComp}(u_1,m)\ \&\ \mathrm{MinComp}(u_2,m) \rightarrow \forall w\ (w\ \varepsilon\ u_1 \leftrightarrow w\ \varepsilon\ u_2)$.

From $M \vDash \mathrm{MinComp}(u,y*a)\ \&\ \mathrm{MinComp}(u',y*a)\ \&\ w'\ \varepsilon\ u'$ , $M \vDash w'\ \varepsilon\ u$.

$\Rightarrow$ from $M \vDash \mathrm{Comp}(u,y*a)\ \&\ \mathrm{Pair}[y*a,z,w']\ \&\ \mathrm{Pair}[y*a,z',w_1]\ \&\ w_1\ \varepsilon\ u$,

$M \vDash z=z'$,

$\Rightarrow$ $M \vDash F_1(y,v,z)$, as required.

This completes the proof of Claim 8.

<u>Claim 9</u>: $QT^+ \vdash \forall y \in J\ \forall v,z \in I\ (H(y,v) \rightarrow (F_2(y,v,z) \rightarrow H(y*b,z)))$.

<u>Claim 10</u>: $QT^+ \vdash \forall y \in J\ \forall v,z \in I\ (H(y,v)\ \&\ H(y*b,z) \rightarrow F_2(y,v,z))$.

These two claims are proved exactly analogously to Claims 8 and 9.

From the definition of the string form J we have

$M \vDash \forall m \in J\ \exists! y \in I\ \exists u \in I\ \exists w \subseteq_p u\ (\mathrm{MinComp}(u,m)\ \&\ \mathrm{Pair}[m,y,w]\ \&\ w\ \varepsilon\ u)$.

So from the definition of H we have

(i) $QT^+ \vdash \forall m \in J\ \exists! y \in I\ H(m,y)$.

From Claims 5 and 6 we have

(iia) $QT^+ \vdash \forall y \in I\ (H(a,y) \leftrightarrow y=p)$,

and

(iib) $QT^+ \vdash \forall y \in I\ (H(b,y) \leftrightarrow y=q)$.

From Claim 7 and 8 we have

(iiia) $QT^+ \vdash \forall y \in J\ \forall v,z \in I\ (H(y,v) \rightarrow (H(y*a,z) \rightarrow F_1(y,v,z)))$,

and, from Claims 9 and 10, we obtain

(iiib) $QT^+ \vdash \forall y \in J\ \forall v,z \in I\ (H(y,v) \rightarrow (H(y*b,z) \rightarrow F_2(y,v,z)))$.

This concludes the proof of the Theorem.■

5.2(a) For any string form $I \subseteq I_\alpha$ and $I \subseteq I_{Add}$ there is a string form $J \equiv I_{Add\alpha} \subseteq I$ such that

$QT^+ \vdash \forall x,y \in J\ \forall u,v,w\ (A^\#(x,u)\ \&\ A^\#(y,v)\ \&\ AddTally(u,v,w) \rightarrow A^\#(x*y,w))$.

<u>Proof</u>: Let J(y) abbreviate

$I(y)\ \&\ \forall x \in I\ \forall u,v,w\ (A^\#(x,u)\ \&\ A^\#(y,v)\ \&\ AddTally(u,v,w) \rightarrow A^\#(x*y,w))$.

Since I may be assumed to be closed under * and downward closed under ≤, we may assume that I is closed under AddTally, α and β.

We argue that J is a string form.

For y=a, we have that M ⊧ I(a).

Assume M ⊧ $A^{\#}$(x,u) & $A^{\#}$(y,v) & AddTally(u,v,w) where M ⊧ I(x).

Then M ⊧ $A^{\#}$(a,v), whence, by (ii$a^{\alpha}$), M ⊧ v=bb.

From M ⊧ $A^{\#}$(x,u), M ⊧ $Tally_b$(u).

Then M ⊧ AddTally(u,bb,w), and by 3.4(d), M ⊧ AddTally(u,bb,Su).

By single-valuedness of AddTally, M ⊧ w=Su.

On the other hand, by ($i^{\alpha}$), M ⊧ ∃!w' ∈ I $A^{\#}$(x*a,w').

From M ⊧ $A^{\#}$(x,u), by (iii$a^{\alpha}$), M ⊧ w'=u*b, and from M ⊧ $Tally_b$(u),

M ⊧ w'=Su. Hence M ⊧ $A^{\#}$(x*a,Su).

Then from M ⊧ w=Su, M ⊧ $A^{\#}$(x*a,w), as required.

For y=b, again we have M ⊧ I(b).

Assume M ⊧ $A^{\#}$(x,u) & $A^{\#}$(y,v) & AddTally(u,v,w) where M ⊧ I(x).

Then M ⊧ $A^{\#}$(b,v). By (ii$b^{\alpha}$), M ⊧ v=b, so M ⊧ AddTally(u,b,w).

By definition of AddTally, M ⊧ w=u.

By ($i^{\alpha}$), M ⊧ ∃!w' ∈ I $A^{\#}$(x*b,w').

Hence from M ⊧ $A^{\#}$(x,u), by (iii$b^{\alpha}$), M ⊧ w'=u. Thus M ⊧ $A^{\#}$(x*b,u).

But then from M ⊧ w=u, M ⊧ $A^{\#}$(x*b,w), as required.

Suppose now that M ⊧ J(y).

Then M ⊧ I(y), whence M ⊧ I(y*a) because I is a string concept.

Assume now that M ⊧ $A^{\#}$(x,u) & $A^{\#}$(y*a,v) & AddTally(u,v,w)

where M ⊧ I(x).

Then M ⊧ $Tally_b$(u) . By ($i^{\alpha}$), M ⊧ ∃!$v_0$ ∈ I ($Tally_b$($v_0$) & $A^{\#}$(y,$v_0$)).

From M ⊨ $A^{\#}(y*a,v)$, by (iii$a^{\alpha}$), M ⊨ $v=v_0*b$.

By 3.5(a), M ⊨ $\exists! w_0 \in I$ $\text{AddTally}_b(u,v_0,w_0)$.

Hence from M ⊨ $A^{\#}(x,u)$ and hypothesis M ⊨ J(y), M ⊨ $A^{\#}(x*y,w_0)$.

From M ⊨ AddTally(u,v,w), M ⊨ AddTally($u,v_0*b$ ,w).

But since M ⊨ $\text{Tally}_b(u)$ & $\text{Tally}_b(v_0)$, from M ⊨ $\text{AddTally}_b(u,v_0,w_0)$, by 3.4(e),

M ⊨ $\text{AddTally}_b(u,v_0*b,w_0*b)$.

Then by single-valuedness of AddTally, M ⊨ $w=w_0*b$.

Since M ⊨ I(y*a), we have from M ⊨ I(x), by ($i^{\alpha}$), that

M ⊨ ∃!w' ∈ I $A^{\#}(x*(y*a),w')$.

But M ⊨ x*(y*a)=(x*y)*a. Hence M ⊨ $A^{\#}((x*y)*a,w')$.

From M ⊨ $A^{\#}(x*y,w_0)$, by (iii$a^{\alpha}$), M ⊨ $w'=w_0*b$, and from M ⊨ $w=w_0*b=w'$,

M ⊨ w=w'.

But then from M ⊨ $A^{\#}(x*(y*a),w')$, M ⊨ $A^{\#}(x*(y*a),w)$, as required.

Therefore, M ⊨ J(y*a).

On the other hand, for yb, we again have, from M ⊨ I(y), that M ⊨ I(y*b).

Assume that M ⊨ $A^{\#}(x,u)$ & $A^{\#}(y*b,v)$ & AddTally(u,v,w) where M ⊨ I(x).

By ($i^{\alpha}$), M ⊨ $\exists! v_0 \in I$ $A^{\#}(y,v_0)$.

Then from M ⊨ $A^{\#}(y*b,v)$, by (iii$b^{\alpha}$), M ⊨ $v=v_0$.

By 3.5(a), M ⊨ $\exists! w_0 \in I$ AddTally($u,v_0,w_0$). So M ⊨ AddTally($u,v,w_0$).

Then from M ⊨ AddTally(u,v,w), by single-valuedness of AddTally,

M ⊨ $w=w_0$.

From hypothesis M ⊨ J(y), M ⊨ $A^{\#}(x*y,w_0)$.

Since $M \models I(y*b)$, we have from $M \models I(x)$, by $(i^{\alpha})$, that

$$M \models \exists! w' \in I \;\; A^{\#}(x*(y*b),w').$$

But $M \models x*(y*b)=(x*y)*b$. So $M \models A^{\#}((x*y)*b,w')$.

But since $M \models I(x*y)$, from $M \models A^{\#}(x*y,w_0)$, by $(iiib^{\alpha})$, $M \models w'=w_0*$.

So from $M \models w'=w_0=w$, $M \models w'=w$.

But then from $M \models A^{\#}(x*(y*b),w')$, $M \models A^{\#}(x*(y*b),w)$, as required.

Therefore, $M \models J(y*b)$, which completes the argument that J is a string form. Then the claim follows immediately.■

6.1(a) $QT^{+} \vdash I^{*}(x) \;\&\; x_2 E x \;\rightarrow\; \forall u,v\, (A^{\#}(x_2,u) \;\&\; B^{\#}(x_2,v) \;\rightarrow\; Sv \leq u)$.

(b) $QT^{+} \vdash I^{*}(x) \;\&\; I^{*}(y) \;\&\; z=bxy \;\rightarrow\; I^{*}(z)$.

(c) $QT^{+} \vdash I^{*}(x) \;\&\; I^{*}(u) \;\&\; bxy=buv \;\rightarrow\; x=u \;\&\; y=v$.

(d) $QT^{+} \vdash I^{*}(x) \;\rightarrow (x \subseteq_p a \leftrightarrow x=a)$.

(e) $QT^{+} \vdash I^{*}(x) \;\&\; I^{*}(y) \;\&\; I^{*}(z) \;\rightarrow\; (x \subseteq_p byz \leftrightarrow x=byz \;\; v \;\; x \subseteq_p y \;\; v \;\; x \subseteq_p z)$.

<u>Proof</u>: (a) Assume $M \models A^{\#}(x_2,u) \;\&\; B^{\#}(x_2,v)$ where $M \models x_2 E x$ and $M \models I^{*}(x)$.

Then $M \models \exists x_1\; x=x_1x_2 \;\&\; x \neq a$, that is, $M \models x_1 B x$.

From $M \models I^{*}(x)$, $M \models J^{*}(x) \;\&\; I^{*}(x_1) \;\&\; Æ(x)$, and also $M \models J^{*}(x_1)$.

By $(i^{\alpha})$ and $(i^{\beta})$, $M \models \exists! v_1 \in J^{*}\; A^{\#}(x_1,v_1) \;\&\; \exists! w_1 \in J^{*}\; B^{\#}(x_1,w_1)$.

Now, from $M \models Æ(x)$, $M \models v_1 \leq w_1$.

Also from $(i^{\alpha})$, $M \models \exists! y \in J^{*}\; A^{\#}(x,y)$, and from $(i^{\beta})$, $M \models \exists! z \in J^{*}\; B^{\#}(x,z)$, and we have that $M \models A^{\#}(x_1*x_2,y) \;\&\; B^{\#}(x_1*x_2,z)$.

On the other hand, since also $M \vDash J^*(x_2)$, again by ($i^\alpha$) and ($i^\beta$) we have

$$M \vDash \exists! v_2 \in J^* A^\#(x_2,v_2) \ \&\ \exists! w_2 \in J^* B^\#(x_2,w_2).$$

By 3.5(a), $M \vDash \exists! z_1 \in J^*(Tally_b(z_1) \ \&\ AddTally(v_1,v_2,z_1))$

and $M \vDash \exists! z_2 \in J^*(Tally_b(z_2) \ \&\ AddTally(w_1,w_2,z_2))$.

Now, from $M \vDash A^\#(x_1,v_1) \ \&\ A^\#(x_2,v_2)$, by 5.2(a),

$$M \vDash A^\#(x_1{*}x_2,z_1),$$

and from $M \vDash B^\#(x_1,w_1) \ \&\ B^\#(x_2,w_2)$, by 5.2(b),

$$M \vDash B^\#(x_1{*}x_2,z_2).$$

On the other hand, from $M \vDash Æ(x) \ \&\ A^\#(x,y) \ \&\ B^\#(x,z)$, $M \vDash y=Sz$.

So from $M \vDash x=x_1x_2$, $M \vDash A^\#(x,z_1) \ \&\ B^\#(x,z_2)$.

Then since $M \vDash J^*(z_1) \ \&\ J^*(z_2)$, $M \vDash z_1=y \ \&\ z_2=z$, and from $M \vDash y=Sz$,

$$M \vDash z_1=Sz_2.$$

Hence from $M \vDash AddTally(v_1,v_2,z_1) \ \&\ AddTally(w_1,w_2,z_2)) \ \&\ v_1 \leq w_1$,

by 3.5(i), $M \vDash Sw_2 \leq v_2$.

From the uniqueness of $v_2,w_2$ we have that $M \vDash Sv \leq u$, as required.

(b) Assume $M \vDash z=bxy$ where $M \vDash I^*(x) \ \&\ I^*(y)$.

Then $M \vDash J^*(x) \ \&\ J^*(y)$, and since $J^*$ is a string form closed under $*$,

$M \vDash J^*(bxy)$.

From $M \vDash I^*(x) \ \&\ I^*(y)$, we have that $M \vDash Æ(x) \ \&\ Æ(y)$. It suffices to show that $M \vDash Æ(z)$.

We proceed to show conditions (c1) and (c2) hold.

By $(i^\alpha)$ and $(i^\beta)$, $M \vDash \exists! v_1 \in J^* A^\#(x,v_1) \& \exists! v_2 \in J^* B^\#(x,v_2)$

and $M \vDash \exists! w_1 \in J^* A^\#(y,w_1) \& \exists! w_2 \in J^* B^\#(y,w_2)$.

From $M \vDash Æ(x) \& Æ(y)$ we have $M \vDash v_1=Sv_2 \& w_1=Sw_2$.

Again by $(i^\alpha)$ and $(i^\beta)$ we have

$$M \vDash \exists! u_1 \in J^* A^\#(z,u_1) \& \exists! u_2 \in J^* B^\#(z,u_2).$$

Then $M \vDash A^\#(b(xy),u_1) \& B^\#(b(xy),u_2)$.

By $(iib^\alpha)$ and $(iib^\beta)$ we have $M \vDash A^\#(b,b) \& B^\#(b,bb)$.

Once again by $(i^\alpha)$ and $(i^\beta)$,

$$M \vDash \exists! v_3 \in J^* A^\#(xy,v_3) \& \exists! v_4 \in J^* B^\#(xy,v_4).$$

By 3.5(a), $M \vDash \exists! p_1 \in J^*(Tally_b(p_1) \& AddTally(b,v_3,p_1))$

and $M \vDash \exists! p_2 \in J^*(Tally_b(p_2) \& AddTally(bb,v_4,p_2))$.

Then from $M \vDash A^\#(b,b)$, by 5.2(a), $M \vDash A^\#(b(xy),p_1)$,

and from $M \vDash B^\#(b,bb)$, by 5.2(b), $M \vDash B^\#(b(xy),p_2)$.

By 3.4(c), $M \vDash AddTally(b,v_3,v_3)$, whence from $M \vDash AddTally(b,v_3,p_1)$, by single-valuedness of Addtally,

$$M \vDash p_1=v_3.$$

By 3.5(c), $M \vDash AddTally(bb,v_4,Sv_4)$, hence from $M \vDash AddTally(bb,v_4,p_2)$, by single-valuedness of Addtally,

$$M \vDash p_2=Sv_4.$$

By 3.5(a), $M \vDash \exists! q_1 \in J^*(Tally_b(q_1) \& AddTally(v_1,w_1,q_1))$

and $M \vDash \exists! q_2 \in J^*(Tally_b(q_2) \& AddTally(v_2,w_2,q_2))$.

Then from $M \vDash A^\#(x,v_1) \& A^\#(y,w_1)$, by 5.2(a),

$M \vDash A^{\#}(xy,q_1)$,

and from $M \vDash B^{\#}(x,v_2)$ & $B^{\#}(y,w_2)$, by 5.2(b),

$M \vDash B^{\#}(xy,p_2)$,

Hence from $M \vDash A^{\#}(xy,v_3)$ & $B^{\#}(xy,v_4)$, by single-valuedness of $A^{\#}$ and $B^{\#}$,

$M \vDash q_1=v_3$ & $q_2=v_4$.

But from $M \vDash A^{\#}(b(xy),u_1)$ & $A^{\#}(b(xy),p_1)$, by single-valuedness of $A^{\#}$,

$M \vDash u_1=p_1$,

and from $M \vDash z=b(xy)$ & $u_1=p_1=v_3=q_1$, $M \vDash A^{\#}(z,q_1)$.

Also, from $M \vDash B^{\#}(b(xy),u_2)$ & $B^{\#}(b(xy),p_2)$, by single-valuedness of $B^{\#}$,

$M \vDash u_2=p_2$.

Then from $M \vDash z=b(xy)$ & $u_2=p_2=Sv_4=Sq_2$, $M \vDash B^{\#}(z,Sq_2)$.

Now, from $M \vDash AddTally(v_1,w_1,q_1)$ & $v_1=Sv_2$ & $w_1=Sw_2$ we have

$M \vDash AddTally(Sv_2,Sw_2,q_1)$.

On the other hand, from $M \vDash AddTally(v_2,w_2,q_2)$, by 3.5(d),

$M \vDash AddTally(Sv_2,w_2,Sq_2)$.

By 3.4(e), $M \vDash AddTally(Sv_2,Sw_2,SSq_2)$.

But then from $M \vDash AddTally(Sv_2,Sw_2,q_1)$ by single-valuedness of Addtally,

$M \vDash q_1=SSq_2$.

Since $M \vDash A^{\#}(z,q_1)$ & $B^{\#}(z,Sq_2)$, this suffices to establish (c1).

For (c2), assume $M \vDash uBz$ & $A^{\#}(u,v_1)$ & $B^{\#}(u,v_2)$.

Then $M \vDash uBb(xy)$, and by 3.7(c),

$$M \vDash u{=}b \ \ v \ \ uBbx \ \ v \ \ u{=}bx \ \ v \ \ \exists y_1(y_1Bu \ \& \ u{=}bxy_1).$$

(1) $M \vDash u{=}b$.

By ($iib^{\alpha}$) and ($iib^{\beta}$), $M \vDash A^{\#}(b,b) \ \& \ B^{\#}(b,bb)$.

Then $M \vDash A^{\#}(u,b) \ \& \ B^{\#}(u,bb)$, and by single-valuedness of $A^{\#}$ and $B^{\#}$,

$$M \vDash v_1{=}b \ \& \ v_2{=}bb.$$

But then $M \vDash v_1{\le}v_2$, as required.

(2) $M \vDash uBbx$.

Then $M \vDash \exists x_1(x_1Bx \ \& \ u{=}bx_1)$.

By ($i^{\alpha}$) and ($i^{\beta}$), $M \vDash \exists! u_1 \in J^* \ A^{\#}(x_1,u_1) \ \& \ \exists! u_2 \in J^* \ B^{\#}(x_1,u_2)$.

From $M \vDash Æ(x)$, $M \vDash u_1{\le}u_2$.

By ($iib^{\alpha}$) and ($iib^{\beta}$), $M \vDash A^{\#}(b,b) \ \& \ B^{\#}(b,bb)$.

By 3.5(a), $M \vDash \exists! p_1 \in J^*(Tally_b(p_1) \ \& \ AddTally(b,u_1,p_1))$

and $M \vDash \exists! p_2 \in J^*(Tally_b(p_2) \ \& \ AddTally(bb,u_2,p_2))$,

Then by 5.2(a), $M \vDash A^{\#}(bx_1,p_1)$, and by 5.2(b), $M \vDash B^{\#}(bx_1,p_2)$.

By 3.4(c), $M \vDash AddTally(b,u_1,u_1)$,

By 3.5(c), $M \vDash AddTally(bb,u_2,Su_2)$.

Hence from $M \vDash AddTally(b,u_1,p_1) \ \& \ AddTally(bb,u_2,p_2)$,

by single-valuedness of Addtally, $M \vDash p_1{=}u_1 \ \& \ p_2{=}Su_2$.

Now, from $M \vDash u{=}bx_1 \ \& \ A^{\#}(u,v_1) \ \& \ B^{\#}(u,v_2)$, we have

$$M \vDash A^{\#}(bx_1,v_1) \ \& \ B^{\#}(bx_1,v_2).$$

Then from $M \vDash A^{\#}(bx_1,p_1) \ \& \ B^{\#}(bx_1,p_2)$, by single-valuedness of $A^{\#}$ and $B^{\#}$,

$$M \vDash v_1{=}p_1 \ \& \ v_2{=}p_2,$$

whence $M \vDash v_1=u_1$ & $v_2=Su_2$.

But then from $M \vDash u_1 \leq u_2$ we have that $M \vDash v_1=u_1<Su_2=v_2$.

By single-valuedness of $A^\#$ and $B^\#$, this suffices to establish (c2) in this case.

(3) $M \vDash u=bx$.

By ($i^\alpha$) and ($i^\beta$), $M \vDash \exists! u_1 \in J^* A^\#(x,u_1)$ & $\exists! u_2 \in J^* B^\#(x,u_2)$.

From $M \vDash Æ(x)$, $M \vDash u_1=Su_2$.

On the other hand, by ($iib^\alpha$) and ($iib^\beta$),

$$M \vDash A^\#(b,b) \ \& \ B^\#(b,bb).$$

By 3.5(a), $M \vDash \exists! p_1 \in J^*(Tally_b(p_1)$ & $AddTally(b,u_1,p_1))$

and $M \vDash \exists! p_2 \in J^*(Tally_b(p_2)$ & $AddTally(bb,u_2,p_2))$.

Reasoning exactly as in (2) with bx in place of $bx_1$ we obtain

$$M \vDash v_1=u_1=Su_2=v_2.$$

By single-valuedness of $A^\#$ and $B^\#$, this suffices.

(4) $M \vDash \exists y_1(y_1By$ & $u=bxy_1)$.

By ($i^\alpha$) and ($i^\beta$), $M \vDash \exists! w_1 \in J^* A^\#(y_1,w_1)$ & $\exists! w_2 \in J^* B^\#(y_1,w_2)$.

From $M \vDash Æ(y)$, $M \vDash w_1 \leq w_2$.

Also by ($i^\alpha$) and ($i^\beta$), $M \vDash \exists! u_1 \in J^* A^\#(x,u_1)$ & $\exists! u_2 \in J^* B^\#(x,u_2)$.

From $M \vDash Æ(x)$, $M \vDash u_1=Su_2$.

By 3.5(a), $M \vDash \exists! q_1 \in J^*(Tally_b(q_1)$ & $AddTally(u_1,w_1,q_1))$

and $M \vDash \exists! q_2 \in J^*(Tally_b(q_2)$ & $AddTally(u_2,w_2,q_2))$,

We then reason as in (1) with u in place of z and $y_1$ in place of y that

$$M \vDash A^\#(u,q_1) \ \& \ B^\#(u,Sq_2).$$

From $M \vDash \mathrm{AddTally}(u_1,w_1,q_1)$ & $u_1=Su_2$,

$$M \vDash \mathrm{AddTally}(Su_2,w_1,q_1).$$

By 3.5(a), $M \vDash \exists! q_3 \in J^*(\mathrm{Tally}_b(q_3)$ & $\mathrm{AddTally}(Su_2,w_2,q_3))$,

whence from $M \vDash w_1 \leq w_2$, by 3.5(b), $M \vDash q_1 \leq q_3$.

From $M \vDash \mathrm{AddTally}(u_2,w_2,q_2))$, by 3.5(d), $M \vDash \mathrm{AddTally}(Su_2,w_2,Sq_2)$.

By single-valuedness of Addtally, $M \vDash q_3=Sq_2$.

Hence $M \vDash q_1 \leq Sq_2$.

Since $M \vDash A^{\#}(u,q_1)$ & $B^{\#}(u,Sq_2)$, this suffices to establish (c2) given single-valuedness of $A^{\#}$ and $B^{\#}$. This completes the argument for $M \vDash Æ(z)$.

(c) Assume $M \vDash bxy=buv$ where $M \vDash I^*(x)$ & $I^*(u)$.

Then $M \vDash J^*(x)$ & $J^*(u)$ & $Æ(x)$ & $Æ(u)$. By (QT3), $M \vDash xy=uv$.

So $M \vDash xB(xy)$ & $uB(xy)$, and by 3.7(a),

$$M \vDash (x=u \ \& \ y=v) \vee xBu \vee uBx.$$

Suppose that

(1) $M \vDash xBu$.

By ($i^{\alpha}$) and ($i^{\beta}$), $M \vDash \exists! x_1 \in J^* A^{\#}(x,x_1)$ & $\exists! x_2 \in J^* B^{\#}(x,x_2)$.

From $M \vDash Æ(u)$, $M \vDash x_1 \leq x_2$.

On the other hand, from $M \vDash Æ(x)$, $M \vDash x_1=Sx_2$.

But then $M \vDash x_1 \leq x_2 < Sx_2 = x_1$, contradicting $M \vDash I_0(x_1)$.

Hence (1) is ruled out.

(2) $M \vDash uBx$.

Ruled out exactly analogously to (a).

Therefore, $M \vDash x=u \,\&\, y=v$, as required.

(d) is immediate from the definition of $\subseteq_p$ by (QT2).■

7.3. WQT* is locally finitely satisfiable.

<u>Proof</u>: Let S be a finite set of axioms of WQT*.

For variable-free terms s, t of $\mathcal{L}_{QT,\sqsubseteq^*}$, let $s \sim t \Leftrightarrow val(s)=val(t)$, that is, if s, t represent the same string. E.g.,

$$a*(b*(a*b)) \sim a*((b*a)*b) \sim (a*b)*(a*b) \sim ((a*b)*a)*b \sim (a*(b*a))*b.$$

$\sim$ being an equivalence relation between terms, we let $[t] = \{ s \mid t \sim s\}$.

Now, let $D = \{a,b,t_1,\ldots,t_n\}$, where $t_1,\ldots,t_n$ are all variable-free terms occurring in S. We let $D^* = \{ [t] \mid t \in D\}$. Since the equivalence classes of terms with respect to $\sim$ can be identified with strings, $D^*$ consists of a, b and the strings represented by terms occurring in S. We take $D^*$ to be the domain of the model, *M*, and let the letters a, b denote [a] and [b], resp. .

Let $f^M : D^* \times D^* \to D^*$, where, for any $[u], [v] \in D^*$,

$f^M([u],[v]) = [t]$ if for some $t \in D$, $t \sim (u*v)$, and $f^M([u],[v]) = b$ otherwise,

interpret the binary operation *.

Let $R^{\mathcal{M}} \subseteq D^* \times D^*$, where, for any $[u], [v] \in D^*$,

$R^{\mathcal{M}}([u],[v]) \Leftrightarrow$ for some $s, t \in \Sigma^\tau$, $u \sim s$ and $v \sim t$ and s is a subterm of t,

interpret the relational symbol $\sqsubseteq^*$.

Suppose now that $[s_1]=[t_1]$ and $[s_2]=[t_2]$, where $[s_1], [s_2], [t_1], t_2] \in D^*$. Then $s_1 \sim t_1$ and $s_2 \sim t_2$, whence $(s_1 * s_2) \sim (t_1 * t_2)$. Suppose further that for some term $s \in D$, $s \sim (s_1 * s_2)$. Then $f^{\mathcal{M}}([s_1], [s_2])=[s]$. But $s \sim (s_1 * s_2) \sim (t_1 * t_2)$. Hence $f^{\mathcal{M}}([t_1], [t_2])=[s]$, and we have $f^{\mathcal{M}}([s_1], [s_2])=f^{\mathcal{M}}([t_1], [t_2])$. Suppose, on the other hand, that for no term $s \in D$, $s \sim (s_1 * s_2)$. Then $f^{\mathcal{M}}([s_1], [s_2])=[b]$. But $(s_1 * s_2) \sim (t_1 * t_2)$, so for no term $s \in D$, $s \sim (t_1 * t_2)$. Hence $f^{\mathcal{M}}([t_1], [t_2])=[b]$, and again $f^{\mathcal{M}}([s_1], [s_2])=f^{\mathcal{M}}([t_1], [t_2])$. Under the same hypothesis $[s_1]=[t_1]$ and $[s_2]=[t_2]$, we have that

$R^{\mathcal{M}}([s_1], [s_2]) \Leftrightarrow$ for some terms $u_1, u_2 \in \Sigma^\tau$, $s_1 \sim u_1$ and $s_2 \sim u_2$ and

$u_1$ is a subterm of $u_2 \Leftrightarrow$ for some terms $u_1, u_2 \in \Sigma^\tau$, $t_1 \sim u_1$ and $t_2 \sim u_2$ and

$u_1$ is a subterm of $u_2 \Leftrightarrow R^{\mathcal{M}}([t_1], [t_2])$.

Thus the definitions of $f^{\mathcal{M}}$ and $R^{\mathcal{M}}$ do not depend on the choice of terms s, t.

A straightforward induction on the complexity of $\mathcal{L}_{\mathrm{C},\sqsubseteq^*}$-terms shows that if t is among the terms in D, then its interpretation $t^{\mathcal{M}}$ is [t]. It is then immediate that the resulting model $\mathcal{M}$ satisfies all of the axioms in the finite set S.■